\documentclass[11pt]{amsart}
\usepackage{amsmath,amsfonts,amssymb,amsxtra,latexsym,amscd,enumerate,amsthm}

\usepackage{latexsym}
\usepackage{epsfig}

\newcommand\ges{\gtrsim}

\newcommand\g{\gamma}
\newcommand\q{\frac{1}{2}}
\newcommand\e{\varepsilon}

\renewcommand\l{\lambda}
\newcommand\D{\Delta}

\newcommand\les{\lesssim}

\newcommand\R{\mathbb{R}}
\newcommand\C{\mathbb{C}}
\newcommand\Z{\mathbb{Z}}
\newcommand\N{\mathbb{N}}
\renewcommand\S{\mathbb{S}}

\newcommand{\lessa}{\lesssim_A}
\newcommand{\inT}{\int_0^T \!\!\!\! \int_{\R^3}\!\!} 

\newcommand{\tU}{\tilde{U}}
\newcommand{\tV}{\tilde{V}}
\newcommand{\tri}{|\!|\!|}

\newtheorem{t1}{Theorem}
\newtheorem{l1}[t1]{Lemma}
\newtheorem{p1}[t1]{Proposition}
\newtheorem{c1}[t1]{Corollary}
\newtheorem{r1}[t1]{Remark}

\newtheorem{d1}[t1]{Definition}

\begin{document}
\title[Maxwell-Schr\"odinger]{Global wellposedness in the energy space
  for the Maxwell-Schr\"odinger system}

\author{Ioan Bejenaru}
\address{ Department of Mathematics \\
  Texas A\&M University}

\author{ Daniel Tataru}
\address {Department of Mathematics \\
  University of California, Berkeley}

\thanks{ The first author was
  partially supported by NSF grant DMS0738442.  \\ The second author was
  partially supported by NSF grant DMS0354539}

\begin{abstract}
  We prove that the Maxwell-Schr\"odinger system in $\R^{3+1}$ is
  globally well-posed in the energy space. The key element of the
  proof is to obtain a short time wave packet parametrix for the
  magnetic Schr\"odinger equation, which leads to linear, bilinear and
  trilinear estimates. These, in turn, are extended to larger time
  scales via a bootstrap argument.

\end{abstract}

\maketitle

\section{Introduction}

The Maxwell-Schr\"odinger system in $\R^{3+1}$ describes the evolution
of a charged non-relativistic quantum mechanical particle interacting
with the classical electro-magnetic field it generates. It has the
form
\begin{equation} \label{ms} \left\{
    \begin{aligned} 
      & iu_{t}-\Delta_A u= \phi u \\
      & - \D \phi + \partial_t \ \text{div} A = \rho, \qquad \qquad   \rho=|u|^2 \\
      & \Box A + \nabla(\partial_t \phi + \text{div} A)= J, 
\qquad J=2 Im (\bar u, \nabla_A u)
    \end{aligned}
  \right.
\end{equation}
where $u$ is the wave function of the particle, $(\phi, A)$ is the
electro-magnetic potential,
\[
(u,A,\phi) : \R^{3} \times \R \rightarrow \C \times \R \times
\R^3
\]
 and $\nabla_A = \nabla - iA$, $\D_A = \nabla_A^2$.

The system is invariant under the gauge transform:
\[
(u',\phi',A') \rightarrow (e^{i\l}u, \phi - \partial_t \l, A+
\nabla\l), \qquad \l : R^3 \times \R \rightarrow \R
\]
where $\l : R^3 \times \R \rightarrow \R$. To remove this degree of
freedom we need to fix the gauge. In this article we choose to work
in the Coulomb gauge
\begin{equation} \label{gaugecond}
\text{div} A =0
\end{equation}
Under this assumption, the system can be rewritten as:
\begin{equation} \label{msc} \left\{
    \begin{aligned} 
      & iu_{t}-\Delta_A u= \phi u \\
      & \Box A =P J
    \end{aligned}
  \right.
\end{equation}
where $\phi=(-\D)^{-1}(|u|^2)$ and $P = 1 - \nabla \text{div} \D^{-1}$
is the projection on the divergence free vectors functions - also
called Helmholtz projection. 
 We consider the above system with a set
of initial data chosen in Sobolev spaces:
\[
(u(0), A(0), A_t(0))= (u_0, A_0, A_1) \in H^s \times H^{\sigma} \times
H^{\sigma-1}
\]
The gauge condition \eqref{gaugecond} is conserved in time provided
the initial data $(A_0, A_1)$ satisfies it due to the form of the
second equation in \eqref{msc}.

The conserved quantities  associated to the system are the charge and
the energy,
\[
Q(u) = \int_{\R^3} |u|^2 dx
\]
\[
 E(u) = \int_{\R^3} |\nabla_A u|^2+
\frac12(|A_t|^2 + |\nabla_x A|^2) + \frac12 |\nabla \phi|^2 dx
\]

The local well-posedness of the system in various Sobolev spaces above
the energy level is known, see \cite{MR1207265}, \cite{MR2181763}.  On
the other hand the existence of weak energy solutions is established
in \cite{MR1331696}.  The main outstanding problem which we seek to
address is the {\em well-posedness in the energy space}. Our result is

\begin{t1}
  The Maxwell-Schr\"odinger system \eqref{msc} is globally well-posed
  in the energy space $H^1 \times H^1 \times L^2$ in the following
  sense:

i) (regular solutions)  For each  initial data 
\[
(u_0, A_0, A_1) \in H^{2} \times H^{2} \times H^{1}
\]
 there exists an unique global solution
\[
 (u,A) \in C(\R, H^2) \times C(\R, H^2) \cap C^1(\R, H^{1}).
\]

ii) (rough solutions) For each  initial data 
\[
(u_0, A_0, A_1) \in H^{1} \times H^{1} \times L^2
\]
 there exists a  global solution
\[
 (u,A) \in C(\R, H^1) \times C(\R, H^1) \cap C^1(\R, L^2).
\]
which is the unique strong limit of the regular solutions in (i).

iii) (continuous dependence) The solutions $(u,A)$ in (ii) depend 
continuously on the initial data in $H^{1} \times H^{1} \times L^2$.
\label{tmain}\end{t1}

We remark that in the process of proving the above results we
establish some additional regularity properties for the energy
solutions $(u,A)$ which suffice both for the uniqueness and the
continuous dependence results. Traditionally these regularity
properties are described using $X^{s,b}$ type spaces. Instead here we
use the related $U^2$ and $V^2$ type spaces associated to both the wave
equation and the magnetic Schr\"odinger equation. These are introduced
in the next section; for more details we refer the reader to
\cite{MR2094851}, \cite{1dnls}, \cite{HHK}.

\begin{r1} We note that in some directions our analysis yields
  stronger results than as stated in the theorem. Precisely, the same
  arguments as those in Section~\ref{stmain} also yield:

a) Local in time a-priori estimates in $H^\beta \times H^1 \times L^2$
for $\beta > \frac12$. This is exactly the range allowed for  $\beta$
in Lemma~\ref{schnonlin} and Lemma~\ref{wnonlin} (a).

b) Local in time well-posedness in $H^\beta \times H^1 \times L^2$ for
$\beta > \frac34$. This reduced range arises due to
Lemma~\ref{wnonlin} (b).
\end{r1}

The nonlinearities on the right hand side of both equations in
\eqref{msc} are fairly mild.  Indeed, if $\Delta_A$ were replaced by
$\Delta$ then it would be quite straightforward to iteratively close
the argument in $X^{s,b}$ or Strichartz spaces.  For the magnetic
potential $A$ it is quite reasonable to hope to obtain an $X^{s,b}$
type regularity. Thus the main difficulty stems from the linear
magnetic Schr\"odinger equation
\begin{equation}
i u_t -\Delta_A u = f, \qquad u(0) = u_0
\label{sca}\end{equation}

The linear and bilinear estimates for $L^2$ solutions to \eqref{sca} are
summarized in Theorem~\ref{ta} in Section~\ref{hss}. The rest of 
the section is devoted to the well-posedness of \eqref{sca} in
$H^2$, $H^{-2}$ and intermediate spaces.

The proof of our main result is completed in the following section.  The
first step is to establish local in time a-priori bounds for solutions
to \eqref{msc}, first in $H^1$ and then in more regular spaces.  This
is done by treating the nonlinearities on the right of the equations
in a perturbative manner. The transition from local in time to global
in time is straightforward, using the conserved energy. The second
step is to establish the continuous dependence on the initial data.
This is a consequence of a  Lipschitz dependence result in a weaker
topology. Precisely, we show that the corresponding linearized equation is 
well-posed in $L^2 \times H^\frac12 \times H^{-\frac12}$.

The rest of the paper is devoted to the study of $L^2$ solutions for
\eqref{sca}, with the aim of proving Theorem~\ref{ta}.  Previous
approaches establish Strichartz estimates with a loss of derivatives
for this equation in a perturbative manner, starting from the free
Schr\"odinger equation.  This no longer suffices for $A$ in the energy
space, and instead one needs to  study directly the dispersive
properties for the linear magnetic Schr\"odinger equation.  
Our approach uses some of the ideas described in \cite{multiscale} and
\cite{1dnls}.

To each dyadic frequency $\lambda$, we associate the time scale
$\lambda^{-1}$. On this time scale we show that at frequency $\lambda$
the equation \eqref{sca} is well approximated by its paradifferential
truncation, which is roughly
\begin{equation}
i u_t -\Delta_{A_{<\sqrt{\lambda}}}  u = f, \qquad u(0) = u_0
\label{sca1}\end{equation}
Following the ideas in \cite{bt}, \cite{mmt0}, in
Sections~\ref{scale1},~\ref{scalel} we obtain a wave packet parametrix
for the equation \eqref{sca1} on the $\lambda^{-1}$ time scale.
This allows us to prove sharp Strichartz and square function
estimates, as well as bilinear $L^2$ bounds and trilinear estimates.

In the last section we extend the linear, bilinear and trilinear
estimates to larger time scales. A brute force summation of the short
time bounds yields unacceptably large constants. Heuristically, the
summation can be improved by taking advantage of the localized energy
estimates for the magnetic Schr\"odinger equation. However these are
not straightforward.  Our idea is to obtain them from a weaker
generalized wave packet decomposition where the localization scales in
position and frequency are relaxed as the time scale is iteratively
increased.

\section{Definitions}

Throughout the paper we use the standard Lebesgue spaces $L^q_x$ and
mixed space-time versions $L^p_t L^q_x$ which are defined in the
standard way. To measure regularity of functions at fixed time we use the standard
Sobolev spaces $H^s_x$.  Additional space-time structures will be
defined in the next section.

We now introduce dyadic multipliers and a Littlewood-Paley
decomposition in frequencies. Throughout the paper the letters
$\lambda,\mu,\nu$ and $\gamma$ will be used to denote dyadic values,
i.e. $\l=2^i$ for some $i \in \N$. 

We say that a function $u$ is localized at frequency $\lambda$ if its
Fourier transform is supported in the annulus $\{|\xi| \in
[\frac{\l}8,8\l]\}$ if $\lambda \geq 2$, respectively in the ball
$\{|\xi| \leq 8\}$ if $\lambda = 1$.
 
By $S_\lambda$ we denote a multiplier  with  smooth symbol
$s_\lambda(\xi)$ which is supported in the annulus $\{|\xi| \in
[\frac{\l}2,2\l]\}$ for $\lambda \geq 2$ respectively in the ball $\{
|\xi| \leq 2\}$ if $\lambda = 1$ and satisfies the bounds
\begin{equation} \label{slm}
|\partial^\alpha s_\l(\xi)| \leq c_\alpha \l^{-|\alpha|}
\end{equation}
By $S_{<\lambda}$ we denote a multiplier with smooth symbol
$s_{<\lambda}(\xi)$ which is supported in the ball $\{ |\xi| \leq
2\lambda\}$, equals $1$ in the ball $\{|\xi| \leq \lambda/2\}$ and
satisfies \eqref{slm}.  All implicit constants in the estimates
involving $S_\lambda$, $S_{<\lambda}$ will depend on finitely many
seminorms of its symbol, i.e. on $c_\alpha$ for $|\alpha| \leq N$, for some
large $N$.

Associated to each $\l$ we also consider $\tilde{S}_\l$ to be a
multiplier whose symbol $\tilde{s}_\l$ satisfies \eqref{slm}, is
supported in $\{ |\xi| \in [\frac{\l}4,4\l]\}$ and equals $1$ in the
support of $\{ |\xi| \in [\frac{\l}2,2\l]\}$. The last condition
implies that
\[
\tilde{S}_\l S_\l = S_\l
\]
Similarly we consider $\tilde{\tilde{S}}_\l$ to be a
multiplier whose symbol satisfies \eqref{slm}, is
supported in $\{ |\xi| \in [{\l}/8,8\l]\}$ and equals $1$ in the
support of $\{ |\xi| \in [{\l}/4,4\l]\}$.

\section{ $V^2$ and $U^2$ type spaces} 
\label{u2v2}

Let $H$ be a Hilbert space. Let $V^{2}H$ be the space of right continuous
$H$ valued functions on $\R$ with bounded 2-variation:
\[
\|u\|^{2}_{V^{2}H} =\sup_{(t_{i}) \in T} \sum_{i}
\|u(t_{i+1})-u(t_{i})\|_H^{2}
\]
where $T$ is the set of finite increasing sequences in $\R$.  

Let $U^{2}H$ be the atomic space defined by the atoms:
\[
u=\sum_{i} h_{i} \chi_{[t_{i},t_{i+1})}, \ \sum_{i} \|h_{i}\|_H^{2}=1
\]
for some $(t_{i}) \in T$. We have the inclusion $U^{2}H \subset V^2H$
but in effect these spaces are very close, and also close to the
homogeneous Sobolev space $\dot H^\frac12$. Precisely, we can bracket
them using homogeneous Besov spaces as follows:
\begin{equation}
  \dot B^{\frac12}_{2,1} \subset U^2 \subset V^2 
\subset \dot B^{\frac12}_{2,\infty}
\label{emb}\end{equation}

We denote by $DU^{2}H$ the space of (distributional) derivatives of
$U^2 H$ functions.  Then there is also a duality relation between
$V^{2}H$ and $U^{2}H$, namely
\begin{equation}
(DU^{2}H)^{*}=V^{2}H
\label{dual1}\end{equation}
For more details on the $U^2$ and $V^2$ spaces we refer
the reader to \cite{MR2094851} and \cite{HHK}.

Given an abstract evolution in $H$, 
\[
iu_t = B(t) u, \qquad u(0) = u_0
\] 
which generates a family of bounded evolution operators 
\[
S(t,s): H \to H, \qquad t,s \in \R
\]
we can define the associated spaces $U^2_B H$, $V^2_B H$
and $DU^2_B H$ by
\[
\| u \|_{V^2_B H} = \| S(0,t) u(t)\|_{V^2 H},
\qquad
\| u \|_{U^2_B H} = \| S(0,t) u(t)\|_{V^2 H} 
\]
respectively
\[
DU^2_B L^2 = \{ (i \partial_t -B) u; \ u \in U^2_B H \}
\]
On occasion we need to compare the above spaces associated 
to closely related operators. For this we use the following

\begin{l1} 
Let $H$ be a Hilbert space and $B(t)$, $C(t)$ two 
families of bounded selfadjoint operators in $H$. Suppose that
for $\e$ sufficiently small we have
\[
\| (B-C)u\|_{DU^2_{B}H} \leq \e \| u\|_{U^2_B H}
\]
Then
\[
\| u\|_{U^2_B H} \approx \|u\|_{U^2_C H}, \qquad \| u\|_{V^2_B H} \approx
\|u\|_{V^2_C H}, \qquad \| f\|_{DU^2_B H} \approx \|f\|_{DU^2_C H}
\]
\label{flowdiff}
\end{l1}
\begin{proof}
By conjugating with respect the $B$ flow
we can assume without any restriction in generality that 
$B=0$. Then we can  solve the equation  
\[
i u_t - C u = 0, \qquad u(0) = u_0
\]
by treating $C$ perturbatively to obtain a solution
\[
u= u_0 + \epsilon u_1(t), \qquad \| u_1\|_{U^2 H} \lesssim \|u_0\|_{H} 
\]
Applying this to each step in $U^2_C$ atoms we obtain 
\[
\|u\|_{U^2 H} \lesssim \|u\|_{U^2_C H}
\]
for arbitrary $u$.

For the converse, applying the above result to each step
in a $U^2 H$ atom we conclude that for each $u \in U^2 H$ 
we can find $u_1 \in U^2 H$ so that
\[
\| u+\epsilon u_1\|_{U^2_C H} + \|u_1\|_{U^2 H} \lesssim \|u\|_{U^2 H}
\]
Iterating this shows that 
\[
\| u\|_{U^2_C H}  \lesssim \|u\|_{U^2 H}
\]
Hence $U^2_C H = U^2 H$.

Consider now $f \in DU^2_C H$. Then $f = iu_t -Cu$ for some $u \in
U^2_C H
= U^2 H$. Since $C$ maps $U^2 H$ to $DU^2 H$ this implies that $f \in
DU^2 H$.
Conversely, if $f \in DU^2 H$ then we can solve the inhomogeneous
equation
\[
i u_t - C u = f, \qquad u(0) = 0
\]
iteratively in $U^2 H$. This gives a solution $u \in U^2 H = U^2_C H$,
therefore $f \in DU^2_C H$. We have proved that $DU^2 H= DU^2_C H$.
Then the  last relation $V^2 H = V^2_C H$ follows by duality.
\end{proof}

Following the above procedure we can associate similar spaces
 to the Schr\"odinger flow by pulling back functions to time
$0$ along the flow, namely
\[
\| u \|_{V^2_\Delta L^2} = \| e^{it\D} u\|_{V^2 L^2},
\qquad
\| u \|_{U^2_\Delta L^2} = \| e^{it\D} u\|_{V^2 L^2} 
\]
The magnetic Schr\"odinger equation has time dependent coefficients,
so we replace the above exponential with the corresponding
evolution operators. We denote by $S^A(t,s)$ the family of evolution
operators corresponding to the equation \eqref{sca}. These are $L^2$
isometries.  Then we define
\[
\| u \|_{V^2_A L^2} = \| S(0,t) u(t)\|_{V^2 L^2},
\qquad
\| u \|_{U^2_A L^2} = \| S(0,t) u(t)\|_{V^2 L^2} 
\]
These spaces turn out to be a good replacement for the
$X^{0,\frac12}$ space associated to the Schr\"odinger equations.
We also define
\[
DU^2_A L^2 = \{ (i \partial_t -\Delta_A) u; \ u \in U^2_A L^2 \}
\]
By \eqref{dual1} we have the duality relation
\[
(DU^2_A L^2 )^* = V^2_A L^2 
\]

 When solving the equation \eqref{sca} we let $f \in DU^2_A L^2$,
and we have the straightforward bound
\begin{equation}
\|u\|_{U^2_A L^2} \lesssim \|u_0\|_{L^2} + \| f\|_{DU^2_A L^2}
\label{u2solve}\end{equation}
In our study of nonlinear equations later on we need to estimate
multilinear expressions in $DU^2_A L^2$. By duality, this is always
turned into multilinear estimates involving $V^2_A L^2$ functions.

Finally, we define similar spaces associated to the wave equation.
The wave equation is second order in time therefore we use a half-wave
decomposition and set
\[
\| u \|_{V^2_{\pm} L^2} = \| e^{\pm it|D|} u\|_{V^2 L^2},
\qquad
\| u \|_{U^2_{\pm} L^2} = \| e^{\pm i t|D|} u\|_{V^2 L^2} 
\]
Then the spaces for the full wave equation are
\[
\| u \|_{U^2_{W} L^2} = \| u\|_{U^2_+ L^2 + U^2_- L^2},
\qquad
\| u \|_{V^2_{W} L^2} = \|  u\|_{V^2_+ L^2 + V^2_- L^2}
\]
For the inhomogeneous term in the wave equation we use the 
space $DU^2_{W} L^2$ with norm
\[
\| f\|_{DU^2_{W} L^2} = \| f\|_{DU^2_+ L^2 \cap  DU^2_- L^2}
\]
Then to solve the inhomogeneous wave equation we use
\[
\| \nabla u\|_{U^2_{W} L^2} \lesssim \|\nabla u(0)\|_{L^2} + \| \Box
u\|_{DU^2_{W} L^2}
\]
Similarly we set 
\[
\| u \|_{U^2_{W} H^s} = \| \langle D_x\rangle^{s}  u\|_{U^2_+ L^2 + U^2_- L^2},
\quad
\| u \|_{V^2_{W} H^s} = \|\langle D_x\rangle^{s}  u\|_{V^2_+ L^2 + V^2_- L^2}
\]
and
\[
\| f\|_{DU^2_{W} H^s} = \| \langle D_x\rangle^{s} f\|_{DU^2 L^2}
\]
Such spaces originate in unpublished work of the second author on the
wave-map equation, and have been successfully used in various contexts
so far, see \cite{MR2094851}, \cite{1dnls},\cite{bt},\cite{HHK}.
The Strichartz estimates for the wave equation turn into embeddings
for $U^2_{W} H^s$ spaces. If the indices $(p,q)$ satisfy
\begin{equation}
\label{pq2}
\frac{1}{p}+\frac{1}{q} = 1,   \qquad 2 < p \leq \infty
\end{equation}
then we have
\begin{equation}
\label{stw}
\| u\|_{L^p L^q} \lesssim \| u\|_{V^2_W H^{\frac{2}p} }
\end{equation}

If we consider frequency localized solutions to the wave equation on a
very small time scale, then the wave equation is
ineffective. Precisely,

\begin{l1} 
Let $B_{<\lambda}$ be a function which is localized at frequency $<
\lambda$. Then
\begin{equation}
\| B_{<\lambda} \|_{U^2_W (I;L^2)} \approx\| B_{<\lambda} \|_{U^2
  (I;L^2)}, \qquad  |I| \leq \lambda^{-1}
\label{pkjh}\end{equation}
\end{l1}

\begin{proof}
This follows from the similar bound for the corresponding half-wave
spaces $U^2_\pm L^2$, and by Lemma~\ref{flowdiff} it is a consequence of the fact that 
for a short time the spatial derivatives in the half wave equation
can be treated perturbatively,
\[
\begin{split}
\| |D_x| B_{<\lambda}\|_{DU^2(I, L^2)} \lesssim &\ \|  |D_x| B_{<\lambda}\|_{L^1(I, L^2)}
\lesssim \l \|  B_{<\lambda}\|_{L^1(I, L^2)}
\\ \lesssim &\  |I| \l\ \|   B_{<\lambda}\|_{L^\infty(I, L^2)}  \lesssim
 |I| \l \ \| B_{<\lambda}\|_{U^2(I, L^2)} 
\end{split}
\]
\end{proof}

The finite speed of propagation for the wave equation allows us to
spatially localize functions in the $U^2_W L^2$ spaces.  For $R > 0$
we consider a covering $(Q_i^R)_{i \in \Z^3}$ of $\R^3$ with cubes of
size $R$. Let $\chi_i^R$ be an associated smooth partition of unity.
Then we have the following result:

\begin{l1} Assume that $I$ is a time interval with $|I| \leq R$. Then: 
\begin{equation} \label{wcloc}
\sum_{i\in \Z^3} \| \chi_{i}^R u \|_{U^2_W (I; L^2)}^2 \les \| u \|^2_{U^2_W (I;L^2)}
\end{equation}
\end{l1}
\begin{proof} 
  By rescaling we can take $R=1$. Without any restriction in
  generality we can also assume that $|I| = 1$.  We prove that the
  result holds for one of the two half-wave spaces, say $U^2_+ (I;
  L^2)$. It is enough to verify \eqref{wcloc} for atoms, and further
  for each step in an atom. Hence we can assume that $u$ 
solves the half wave equation
\[
(i \partial_t +|D_x|) u = 0
\]
Then we have 
\[
(i \partial_t +|D_x|) (\chi_{i}^R u ) = [ |D_x|,\chi_{i}^R] u
\]
By standard commutator estimates we have at fixed time
\[
\sum_{i\in \Z^3} \| [ |D_x|,\chi_{i}^R] u(t) \|_{L^2}^2 \lesssim
\|u(t)\|_{L^2}^2
\]
Then  
\[
 \sum_{i\in \Z^3} \| \chi_{i}^R u \|_{U^2_W (I; L^2)}^2
\lesssim \sum_{i\in \Z^3} \| \chi_{i}^R u(0) \|_{ L^2}^2 + \| [
|D_x|,\chi_{i}^R] u(t) \|_{L^1 L^2}^2 \lesssim \|u(0)\|_{L^2}^2
\] 
\end{proof}

Due to the atomic structure, in many estimates it is convenient to
work with the $U^2$ type spaces instead of $V^2$. In order to transfer
the estimates from $U^2$ to $V^2$ we use the following result from
\cite{HHK}:

\begin{p1} 
\label{phhk}
Let $2 < p < \infty$.  If $u \in V^2 H$ then for
  each $0 < \e < 1$ there exist $u_1 \in U^2 H$ and $u_2 \in U^p H$
  such that $u=u_1+u_2$ and
\begin{equation} 
\label{v2u2} |\ln{\e}|^{-1} \| u_1 \|_{U^2 H} + \e^{-1} \|
  u_2 \|_{U^p H} \les \| u \|_{V^2 H} 
\end{equation}
\end{p1}
Here $U^p$ is defined in the same manner as $U^2$ but with the $l^2$
summation replaced by an $l^p$ summation.  One way we use this result 
is as follows:

\begin{c1} 
Let $\e > 0$ and $N$ arbitrarily large. Then 
\[
V^2_W H^s \subset U^2_W H^{s-\epsilon} + L^\infty H^N
\]
\label{chhk}\end{c1}

Following is another example of how this result can be applied.
Typically in our analysis we prove dyadic trilinear estimates of the
form
\begin{equation}
 \label{uuu} \left| \int_I \int_\R^3  S_{\l_1} u  S_{\l_2} \bar{v}, S_{\l_3} B dx
  dt \right| \les C_1(|I|,\lambda_{123})
 \| u \|_{U^2_A L^2} \| v \|_{U^2_A L^2} \| B
\|_{U^2_W L^2}
\end{equation} 
What we need instead is an estimate with  one 
$U^2$ replaced by a $V^2$, say
\begin{equation} \label{uuv} 
 \left| \int_I \int_\R^3  \! S_{\l_1} u  S_{\l_2} \bar{v}, S_{\l_3} B dx
  dt \right| \les C_2(|I|,\lambda_{123})
 \| u \|_{U^2_A L^2} \| v \|_{U^2_A L^2} \| B
\|_{V^2_W L^2}
\end{equation}
Denoting 
\[
\lambda_{max} = \max\{ \lambda_1,\lambda_2,\lambda_3\},
\]
due to Proposition~\ref{phhk} we can easily show that
\begin{l1} \label{principle} Assume \eqref{uuu} holds and $|I| \leq 1$.
Then \eqref{uuv} holds with
\[
C_2(|I|,\lambda_{123}) = C_1(|I|,\lambda_{123}) \ln{\lambda_{max}}.
\]
The same holds if the $V^2$ structure is placed on any of the other two
factors in \eqref{uuv}.
\end{l1}

\begin{proof}
  Without any restriction in generality we assume that $\lambda_1$,
  $\lambda_2$ and $\lambda_3$ are so that the integral in \eqref{uuu}
  is nontrivial. Taking $u$, $v$ and $B$ to be time independent
  frequency localized bump functions we easily see that
\begin{equation}
C_1(|I|,\lambda_{123}) \gtrsim |I| \lambda_{max}^{-N}
\label{pre}\end{equation}
for some sufficiently large $N$.

For each $0 < \epsilon \leq 1$ we decompose $B=B_1+B_2$
as in Proposition~\ref{phhk}.
For $B_1$ we use \eqref{uuu} while for $B_2$ we use 
Bernstein's inequality to estimate 
\[
\begin{split}
\left| \int_I \int_\R^3  \! S_{\l_1} u  S_{\l_2} \bar{v}, S_{\l_3} B_2 dx
  dt \right| \les &\ |I| \lambda_{max}^N
 \| u \|_{L^\infty L^2} \| v \|_{L^\infty L^2} \| B_2
\|_{L^\infty L^2}
\\ \les &\ |I| \lambda_{max}^N
 \| u \|_{U^2_A L^2} \| v \|_{U^2_A L^2} \| B_2
\|_{U^p_W L^2}
\end{split}
\]
Adding the $B_1$ and the $B_2$ bounds gives
\[
C_2(|I|,\lambda_{123}) \lesssim |\ln \epsilon| C_1(|I|,\lambda_{123})
+ \epsilon  |I| \lambda_{max}^N
\]
We set $\epsilon = \lambda_{max}^{-2N}$. Then the conclusion of the
Lemma follows due to \eqref{pre}.

\end{proof}

\section{ The linear magnetic Schr\"odinger equation}
\label{hss}

In this section we summarize the key properties of 
solutions to the homogeneous and inhomogeneous 
linear magnetic Schr\"odinger equation
\begin{equation}
i u_t -\Delta_A u = 0, \qquad u(0) = u_0
\label{schra}\end{equation}
\begin{equation}
i u_t -\Delta_A u = f, \qquad u(0) = u_0
\label{schrainhom}\end{equation}
in $L^2$, and we use them in order to show that the 
above equation is also well-posed in $H^2$, $H^{-2}$ 
and in intermediate spaces. 

We assume that $A \in U^2_W H^1$ with $ \nabla \cdot A = 0$.  All
constants in the estimates depend on the $U^2_W H^1$ norm of $A$ which
is why we introduce the notation $X \lessa Y $, which means $X \leq
C(\|A\|_{U^2_W H^1})Y$.  

The trilinear estimates are concerned with integrals of the form
\[
I^T_{\lambda_1,\lambda_2,\lambda_3}(u,v,B) = \int_0^T \int_{\R^3}
S_{\lambda_1} u \, S_{\lambda_2} \bar v \, S_{\lambda_3} B \ dx dt
\]
where $u$ and $v$ are associated to the magnetic Schr\"odinger
equation and $B$ is associated to the wave equation. In order for the
above integral to be nontrivial the two highest frequencies need to be
comparable. Thus by a slight abuse of notation in the sequel we 
restrict ourselves to the case
\[
\{ \lambda_1, \lambda_2,\lambda_3\} = \{ \lambda, \lambda,\mu\},
 \qquad \mu \leq \lambda
\]
With these notations we have

\begin{t1} \label{ta} For each $A \in U^2_W H^1$ with $ \nabla \cdot A
  = 0$ the equation \eqref{schra} is well-posed in $L^2$.  For each
  $\epsilon > 0$ there exists $\delta > 0$ so that  the following
  properties hold for $0 < T \leq 1$:

(i) Strichartz estimates:
\begin{equation}
\| S_\lambda u\|_{L^p(0,T; L^q)} \lessa T^{\frac{\delta}p} \lambda^{\frac1p} \|
u\|_{U^2_A L^2}
\qquad \frac{2}{p} + \frac{3}q = \frac32,\quad  2 \leq p \leq \infty
\label{seu2}\end{equation} 

(ii) Local energy estimates. For any spatial cube $Q$ of size  $1$ we have 
\begin{equation}
\| S_\lambda u\|_{L^2(0,T; Q)} \lessa T^\delta \lambda^{-\frac12+\e} \|
u\|_{U^2_A L^2}.
\label{le}\end{equation} 

(iii) Local Strichartz estimates. For any spatial cube $Q$ of size $1$ we have:
\begin{equation}
\| S_\lambda u\|_{L^2(0,T; L^6(Q))} \lessa  T^\delta \lambda^{\epsilon} \|
u\|_{U^2_A L^2}.
\label{seu2loc}\end{equation} 

(iv) Trilinear estimates. For any $0 < T \leq 1$  and $\mu \leq
\lambda$ we have
\begin{equation}
|I^T_{\lambda,\lambda,\mu}(u,v,B)|
\lessa  T^\delta \l^\epsilon \min{(1,\mu \l^{-\frac12})} \|u\|_{U^2_A L^2}  \|v\|_{U^2_A L^2}
\|B\|_{U^2_W L^2}.   
\label{llm}\end{equation}
On the other hand if $\mu \ll \lambda$ then
\begin{equation}
|I^T_{\lambda,\mu,\lambda}(u,v,B)| \lessa  T^\delta
\lambda^{-\frac12+\epsilon} \mu^\frac12 \|u\|_{U^2_A L^2}  \|v\|_{U^2_A L^2}
\|B\|_{U^2_W L^2}  
\label{lml}\end{equation}
\end{t1}

The proof of this theorem is quite involved and is relegated to
Sections~\ref{scale1}-\ref{largescale}.

The smallness given by the $T^\delta$ factor is needed in several 
proofs which use either the contraction principle or bootstrap
arguments. However, this factor is nontrivial only in \eqref{seu2}.
Indeed, we have

\begin{r1}
Assume that the conclusion of Theorem~\ref{ta}  holds without the 
$T^\delta$ factor in \eqref{le}, \eqref{seu2loc}, \eqref{llm} and
\eqref{lml}. Then   the conclusion of Theorem~\ref{ta}  holds 
in full.
\label{tdelta}\end{r1}
\begin{proof}
For \eqref{le} we observe that
\[
\| S_\lambda u\|_{L^2(0,T; Q)} \lesssim T^\frac12 \|u\|_{L^\infty L^2}
\lesssim  \| u\|_{U^2_A L^2}.
\]
Interpolating this with \eqref{le} without the $T^\delta$ factor
yields \eqref{le} with a $T^\delta$ factor. By Bernstein's inequality
the same argument works for \eqref{seu2loc}.

For \eqref{llm} we  can also write the obvious estimate
\[
|I^T_{\lambda,\lambda,\mu}(u,v,B)|
\lessa  T \mu^\frac32  \|u\|_{L^\infty L^2}  \|v\|_{L^\infty L^2}
\|B\|_{L^\infty L^2}.   
\]
which is then interpolated with \eqref{llm} without the $T^\delta$
factor. The same argument applies for \eqref{lml}.
\end{proof}

Next we turn our attention to the $H^2$ and $H^{-2}$ well-posedness
for \eqref{schra}.  We make the transition from $L^2$ to $H^2$ and
$H^{-2}$ using the coercive elliptic operator $1 - \Delta_A$. Its properties
are summarized in the following

\begin{l1} \label{1-d}
  For each $0 \leq s \leq 2$ the operator $1 - \Delta_A$ is a
  diffeomorphism
\[
1 - \Delta_A: H^s \to H^{s-2}
\]
which depends continuously on $A \in H^1$.
\end{l1}
The proof uses standard elliptic arguments and is left for the reader.
\bigskip

Using the above operator we define the spaces $U^2_A H^2$, $V^2_A
H^2$, respectively $DU^2_A H^2$ by
\[
\| u\|_{\tU^2_A H^2} = \|(1 -\Delta_A) u\|_{U^2_A L^2},
\qquad
\| u\|_{\tV^2_A H^2} = \|(1 -\Delta_A)u\|_{V^2_A L^2}  
\]
respectively 
\[
\| f\|_{D\tU^2_A H^2} = \|(1-\Delta_A) f\|_{DU^2_A L^2} 
\]

\begin{r1}
  The reason we use the $\tU$, $\tV$ notation above is to differentiate these
  spaces from the $U^2_A H^2$, $V^2_A H^2$, $DU^2_A H^2$ spaces which
  should be defined as in the previous section, with respect to the
  $H^2$ flow of \eqref{schra}. This is not possible at this point, as
  we have not yet proved that \eqref{schra} is well-posed in $H^2$.
  However, after we do so we will prove that the above two sets of
  norms are equivalent. After that the $\tU$, $\tV$ notation is dropped.
\end{r1}

We can transfer the estimates from Theorem~\ref{ta} to the $U^2_A H^2$
spaces by making an elliptic transition between $\Delta_A$ and
$\Delta$:
 
\begin{l1} \label{lh2a}
Let $A \in  U^2_W H^1$ with $ \nabla \cdot A = 0$. Then for each $\e >
0$ there exists $\delta > 0$ so that the following properties hold:

(i) For $p,q$ as in \eqref{seu2} we have the Strichartz estimate
\begin{equation} \label{h2str}
\| S_\lambda u \|_{L^p(0,T; L^q)} \lessa T^{\frac{\delta}p}
\lambda^{-2+\frac1p}  \|u\|_{\tU^2_A H^2}
\end{equation}

(ii) Elliptic representation. Each $u \in U^2_A H^2$ can be
expressed as
\begin{equation}
u = (1-\Delta)^{-1}  (u_e + u_r), \qquad \|u_e\|_{U^2_A L^2} +
T^{-\delta}\|u_r\|_{L^2(0,T; H^{1-\epsilon})} \lessa \|u\|_{\tU^2_A H^2}
\label{ereph2}\end{equation}

(iii) Local energy estimates. 
For any spatial cube $Q$ of size  $1$ we have:
\begin{equation}
\| S_\lambda u\|_{L^2([0,T]\times Q)} \lessa  T^\delta \lambda^{-2 - \frac12 +\epsilon} 
\| u_0\|_{\tU^2_A H^2}
\label{leh2}\end{equation} 

(iv) Local Strichartz estimates. For any spatial cube $Q$ of size  $1$ we have:
\begin{equation}
\| S_\lambda u\|_{L^2(0,T; L^6(Q))} \lessa T^\delta  \lambda^{-2+\epsilon} 
\| u_0\|_{\tU^2_A H^2}
\label{seuloch2}\end{equation}  
\end{l1}

\begin{proof}
(i) The Strichartz estimate \eqref{h2str} follows from
\begin{equation} \label{helpsrt}
\| S_\lambda (1-\Delta) u \|_{L^p L^q} \lessa T^{\frac{\delta}p}
\lambda^{\frac1p}  \|u\|_{\tU^2_A H^2}
\end{equation}
We use the identity
\begin{equation}
(1-\Delta) u = (1- \Delta_A) u - 2i A \nabla u - A^2 u = (1- \Delta_A) u -R_A(u)
\label{deltau}\end{equation}
and estimate each of the three terms.  From the definition of $\tU^2_A
H^2$ and \eqref{seu2} we have:
\[
\| S_\l (1-\Delta_A) u\|_{L^p L^q} \lessa T^{\frac{\delta}p}\l^{\frac1{p}} \| u \|_{\tU^2_A H^2} 
\]
For the second term we use Bernstein's inequality 
and the exponent  relation in \eqref{seu2} to estimate
\[
\begin{split}
\| S_\l (A \nabla u) \|_{L^p L^q} &\ 
\les_A  T^\frac1p \l^{\frac2{p}} \| S_\l (A \nabla u)
\|_{L^\infty L^2} 
\\
&\ \les  T^\frac1p \l^{\frac1{p}} \| A \nabla u \|_{L^\infty H^\frac12}
\\
&\ \les  T^\frac1p \l^{\frac1{p}} \| A\|_{L^\infty H^1} 
\|  u \|_{L^\infty H^2}
\end{split}
\]
Similarly for the last term we obtain
\[
\| S_\l (A^2 u) \|_{L^p L^q} \les  T^\frac1p \l^{\frac1{p}} \| A \|^2_{L^{\infty} H^1} \| u \|_{L^\infty H^2}
\]
This concludes the proof for \eqref{helpsrt} which implies \eqref{h2str}. 

(ii) By \eqref{deltau} we can set
\[
u_e = (1 - \Delta_A) u, \qquad u_r = -R_A(u) = -2iA \nabla u - A^2 u
\]
Hence it remains to prove that
\begin{equation}
\| S_\lambda (A \nabla u) \|_{L^2} + \| S_\lambda (A^2 u)\|_{L^2}
\lessa  T^\delta \lambda^{-1+\epsilon} \|u\|_{\tU^2_A H^2}
\label{rest}\end{equation}
By the argument in Remark~\ref{tdelta}, here and for the rest of the
proof of the lemma we can neglect the $T^\delta$ factors.

We decompose the expression $S_\lambda (A \nabla u)$ as
 \begin{equation} \label{dec100} 
S_\l (A \nabla u) =  S_\l \left(\sum_{\g \les \l} S_\gamma A
 S_\l   \nabla u + \sum_{\g \les \l} S_\l A S_\gamma \nabla u +
   \sum_{\g \ges \l} S_\g A S_\g  \nabla u \right)
\end{equation}
For exponents $(p,q)$ satisfying \eqref{pq2} we invoke the Strichartz
estimates \eqref{stw} for the wave equation.  By Bernstein's
inequality and \eqref{h2str} we can derive a similar bound for $u$,
\begin{equation}
\| S_\gamma \nabla u \|_{L^q L^p} \lesssim_A 
\gamma^{-\frac{2}p} \|u\|_{\tU^2_A H^2}.
\label{sts}\end{equation}
The $L^2$ bound for the product is obtained by multiplying the last two 
inequalities. For the first term in \eqref{dec100} we take $q$ close
to $\infty$, for the second we take $p=\infty$, while for the third
any choice will do.

For the expression $S_{\lambda} (A^2 u)$ we take a triple
Littlewood-Paley decomposition,
\[
S_{\lambda} (A^2 u)= \sum_{\lambda_1,\lambda_2,\lambda_3}
S_{\lambda} (A_{\lambda_1} A_{\lambda_2} u_{\lambda_3})
\]
Then we must have $\lambda_{max} \geq \lambda$. 

If $\lambda_{max} = \lambda_3$ then we use the above Strichartz
inequalities and Bernstein's inequality to estimate the triple product
as
\[
\begin{split}
  \| A_{\lambda_1} A_{\lambda_2} u_{\lambda_3}\|_{L^2} &\ \lesssim
  \|A_{\lambda_1} \|_{L^4 L^\infty} \|A_{\lambda_2} \|_{L^4
    L^\infty} \|u_{\lambda_3}\|_{L^\infty L^2} \\ &\ \lessa 
  \lambda_1^{\frac14} \lambda_2^{\frac14} \lambda_3^{-2}
  \|A_{\lambda_1} \|_{U^2_W H^1} \|A_{\lambda_2} \|_{U^2_W H^1}
  \|u_{\lambda_3}\|_{\tU^2_A H^2}
\end{split}
\]
The summation with respect to $\lambda_1,\lambda_2,\lambda_3$ is
straightforward.

If $\lambda_3 \ll \lambda_{max}$ then there are two possibilities.
One is $\lambda_{max} = \lambda$, in which case we assume w.a.r.g.
that $\lambda = \lambda_1 \geq \lambda_2,\lambda_3$ and estimate
 \[
\begin{split}
\| A_{\lambda_1} A_{\lambda_2} u_{\lambda_3}\|_{L^2} &\ \lesssim
\|A_{\lambda} \|_{L^p L^q} \|A_{\lambda_2} \|_{L^q L^p} 
\|u_{\lambda_3}\|_{L^\infty}
\\ &\ \lessa \lambda_1^{-\frac2p} \lambda_2^{-\frac2q} \lambda_3^{-\frac12}
    \|A_{\lambda_1} \|_{U^2_W H^1} \|A_{\lambda_2} \|_{U^2_W H^1} 
\|u_{\lambda_3}\|_{\tU^2_A H^2}
\end{split}
\]
with $q$ close to infinity. The other possibility is $\lambda_{max}
\gg \lambda$, in which case we must have $\lambda_1 = \lambda_2 \gg
\lambda_3$. Then we estimate the triple product as above, but the
choice of $p$ and $q$ is no longer important. The proof of
\eqref{rest} is concluded.

(iii) We use the representation in \eqref{ereph2}. The bound for $u_r$ holds
without any localization. For $u_e$ we can write
\[
S_\lambda  (1-\Delta)^{-1} u_e = \lambda^{-2}  ( \lambda^2 (1-\Delta)^{-1} \tilde S_\lambda)
S_\lambda u_e
\]
where the symbol of $ \tilde S_\lambda$ is still supported at
frequency $\lambda$ but equals $1$ in the support of the symbol of
$S_\lambda$. The operator $ \lambda^2 \Delta^{-1} \tilde S_\lambda$
is a unit mollifier acting on the $\lambda^{-1}$ scale, therefore 
it is bounded in $l^\infty_Q L^2([0,1] \times Q)$.

(iv) We use the representation in \eqref{ereph2}. By Bernstein's
inequality the bound for $u_r$ holds without any localization.  For
$u_e$ we argue as above.
\end{proof}

Next we define similar spaces $\tU^2_A H^{-2}$, $\tV^2_A H^{-2}$,
respectively $D\tU^2_A H^{-2}$ in a manner similar to the $H^2$ case,
namely
\[
\| u\|_{\tU^2_A H^{-2}} =  \| (1-\Delta_A)^{-1} u \|_{U^2_A L^2}, \qquad
| u\|_{\tV^2_A H^{-2}} =  \| (1-\Delta_A)^{-1} u \|_{V^2_A L^2},
\]
 respectively 
\[
\| f\|_{D\tU^2_A H^{-2}} = \| (1-\Delta_A)^{-1} f \|_{DU^2_A L^2} 
\]
Due to the duality relation~\ref{dual1} we have the $H^2 - H^{-2}$ duality
\begin{equation}
(D\tU^2_A H^{-2})^* = \tV^2_A H^2, 
\qquad (D\tU^2_A H^{2})^* = \tV^2_A H^{-2}
\label{dualh2}\end{equation}
For functions in $U^2_A H^{-2}$ we can prove results similar to
Lemma~\ref{lh2a}:

\begin{l1} \label{lh-2a}
Let $A \in  U^2_W H^1$ with $ \nabla \cdot A = 0$. Then for each $\e >
0$ there exists $\delta > 0$ so that the following properties hold:

(i) For $p,q$ as in \eqref{seu2} we have the Strichartz estimate
\begin{equation} \label{h-2str}
\| S_\lambda u \|_{L^p(0,T; L^q)} \lessa 
T^\frac{\delta}p \lambda^{2+\frac1p}  \|u\|_{\tU^2_A H^{-2}}
\end{equation}

(ii) Elliptic representation. Each $u \in \tU^2_A H^2$ can be
expressed as
\begin{equation}
u = (1-\Delta) (u_e + u_r), \qquad \|u_e\|_{U^2_A L^2} +
T^{-\delta}\|u_r\|_{L^2 (0,T;H^{1-\epsilon})} \lessa \|u\|_{\tU^2_A H^{-2}}
\label{ereph-2}\end{equation}

(iii) Local energy estimates. 
For any spatial cube $Q$ of size  $1$ we have:
\begin{equation}
\| S_\lambda u\|_{L^2([0,T]\times Q)} \lessa  T^\delta \lambda^{2 - \frac12 +\epsilon} 
\| u_0\|_{\tU^2_A H^{-2}}
\label{leh-2}\end{equation} 

(iv) Local Strichartz estimates. For any spatial cube $Q$ of size  $1$ we have:
\begin{equation}
\| S_\lambda u\|_{L^2(0,T; L^6(Q))} \lessa  T^\delta \lambda^{2+\epsilon} 
\| u_0\|_{\tU^2_A H^{-2}}
\label{seuloch-2}\end{equation}  
\end{l1}

\begin{proof} The proof  is similar to the proof of Lemma~\ref{lh2a},
  so we merely outline it. We denote $v = (1-\Delta_A)^{-1} u$. Then
\[
u = (1-\Delta) v - R_A(v)
\]
To prove \eqref{h-2str} we use \eqref{seu2} for $(1-\Delta) v$ 
and it remains to show that
\[
\| S_\lambda \nabla (A v)\|_{L^p L^q}     +\|S_\lambda(A^2
v)\|_{L^p L^q} \lessa T^\delta
\lambda^{2+\frac1p}  \|v\|_{U^2_A L^2}
\]
which is obtained using only the energy estimates for $A$ and $v$.

For \eqref{ereph-2} we set
\[
u_e = v, \qquad u_r = (1-\Delta)^{-1} R_A(v)
\]
Then it remains to show that
\[
\| S_\lambda \nabla (A v)\|_{L^2}     +\|S_\lambda(A^2 v)\|_{L^2} \lessa 
T^\delta \lambda^{1+\e}  \|v\|_{U^2_A L^2}
\]
But this is obtained in the same manner as \eqref{rest} from the
Strichartz estimates for $A$ and $u_2$.

Finally, \eqref{leh-2} and \eqref{seuloch-2} are proved exactly as in
the previous lemma.
\end{proof}

Next we turn our attention to the trilinear bounds, namely the $H^2$
and $H^{-2}$ counterparts of  \eqref{llm} and \eqref{lml}. For
uniformity in notations we set $\tU^2_A L^2 = U^2_A L^2$. Then

\begin{l1} \label{trih2}
Let $k,l \in \{ -2,0,2\}$. Then  for any $0 < T \leq 1$  and $\mu \leq
\lambda$ we have
\begin{equation}
\!\! |I^T_{\lambda,\lambda,\mu}(u,v,B)|
\! \lessa\! T^\delta  \l^\epsilon \min{(1,\mu \l^{-\frac12})} \l^{-k-l} \|u\|_{\tU^2_A H^k}  \|v\|_{\tU^2_A H^l}
\|B\|_{U^2_W L^2}  
\label{llmh2h2}\end{equation}
while if $\mu \ll \lambda$ then
\begin{equation}
|I^T_{\lambda,\mu,\lambda}(u,v,B)| \lessa  T^\delta 
\lambda^{-\frac12+\epsilon} \mu^\frac12 \mu^{-l} \l^{-k} \|u\|_{\tU^2_A H^k}  \|v\|_{\tU^2_A H^l}
\|B\|_{U^2_W L^2}  
\label{lmlh2h2}\end{equation}
\end{l1}

\begin{proof}
The $T^\delta$ factor can be neglected by the argument in Remark~\ref{tdelta}.
We represent 
\[
u = (1-\Delta)^{-\frac{k}2}(u_e + u_r)
\]
where $u_e$, $u_r$ are chosen as in \eqref{ereph2} if $k=2$, as in
\eqref{ereph-2} if $k = -2$ and with $u_r = 0$ if $k=0$. Similarly we set
\[
v = (1-\Delta)^{-\frac{k}2}(v_e + v_r)
\]

We begin with \eqref{llmh2h2} and consider all four combinations.
The estimate for $u_e$ and $v_e$ is exactly \eqref{llm}. The estimate 
for $u_e$ and $v_r$ reads
\[
 |I^1_{\lambda,\lambda,\mu}(u_e,v_r,B)|
 \lessa  \l^\epsilon \min{(1,\mu \l^{-\frac12})}
\|u_e\|_{U^2_A L^2}  \|v_r\|_{L^2 H^{1-\epsilon}}
\|B\|_{U^2_W L^2}  
\]
and is a consequence of the stronger bilinear $L^2$ estimate
\begin{equation}
\| S_\lambda u_e S_\mu B\|_{L^2}  \lessa  \mu^{\frac12} \l^{\epsilon}
\|u_e\|_{U^2_A L^2} 
\|B\|_{U^2_W L^2}  
\label{ulbm}\end{equation}
Due to the finite speed of propagation for the wave equation, see \eqref{wcloc}, 
we can localize this to the unit spatial scale. But on a unit cube $Q$
we use the local Strichartz estimate \eqref{seu2loc} for $u_e$
and the energy estimate for $B$.

The estimate for $u_r$ and $v_e$ is similar. The estimate for 
$u_r$ and $v_r$ reads
\[
|I^1_{\lambda,\lambda,\mu}(u_r,v_r,B)|
 \lessa  \l^\epsilon \min{(1,\mu \l^{-\frac12})}
  \|u_r\|_{L^2 H^{1-\epsilon}}  \|v_r\|_{L^2 H^{1-\epsilon}}
\|B\|_{U^2_W L^2}  
\]
and is easily proved using $L^2$ bounds for $S_\l u_r$, $S_\l v_r$
and an $L^\infty$ bound for $S_\mu B$.

The proof of \eqref{lmlh2h2} is similar. For later use 
we note the bilinear $L^2$ estimates, namely
\begin{equation}
\| S_\mu v_e S_\lambda B\|_{L^2}  \lessa  \mu^{\frac12+\epsilon}
\|u_e\|_{U^2_A L^2} 
\|B\|_{U^2_W L^2}  
\label{umbl}\end{equation}
respectively
\begin{equation}
\| S_\mu (S_\lambda u_e S_\lambda B)\|_{L^2}  \lessa  \mu^{\frac12+\epsilon}
\|u_e\|_{U^2_A L^2} 
\|B\|_{U^2_W L^2}  
\label{ulbl}\end{equation}
Both are proved by localizing to a unit spatial scale and then by
combining the local Strichartz estimate \eqref{seu2loc} for
$u_e$ and $v_e$  and the energy estimate \eqref{stw} for $B$.

\end{proof}

Now we consider the $H^2$ well-posedness of \eqref{schra}.

\begin{p1}
 Let $A \in U^2_W$ with $ \nabla \cdot A = 0$. 
Then the equation \eqref{schrainhom} is well-posed in
 $H^2$.  In addition we have
\[
U^2_A H^2 = \tU^2_A H^2, \qquad V^2_A H^2 = \tV^2_A H^2,
\qquad  DU^2_A H^2 = \tU^2_A H^2
\]
with equivalent norms.
\end{p1}
\begin{proof}
  In order to solve the equation \eqref{schra} with initial data $u_0
  \in H^2$ we consider the equation for $v=(1-\Delta_A) u$ which has
  the form
\[
(i \partial_t - \D_A)  v  =  2(A_t \nabla  -i A A_t) u 
\]
Expressing $u$ in terms of $v$ we obtain 
\begin{equation} \label{eqd}
(i \partial_t - \D_A)  v  =  2 ( A_t \nabla -
i A A_t ) ( 1-\Delta_A)^{-1} v 
\end{equation}
We seek to solve this equation perturbatively in $U^2_A L^2$.
For this we need first to establish suitable mapping properties 
for the operator $A_t \nabla  -i A A_t$.

\begin{l1}
 The operator $A_t \nabla - i A A_t$ satisfies the space-time bound 
\begin{equation}
\|  (A_t \nabla  - 2i A A_t) u \|_{DU^2_A L^2} \lessa T^\delta \| u\|_{\tU^2_A H^2}
\label{atbd}\end{equation}

\end{l1}
\begin{proof}
 By duality, \eqref{atbd}  follows from the bounds
\begin{equation}
\left| \int_{0}^T \int_{\R^3} B \nabla u \bar v  dx dt \right| \lessa T^\delta
 \| u\|_{\tU^2_A H^2}  \| v\|_{V^2_A L^2}\|  B\|_{U^2_W L^2}
\label{h2a} \end{equation}
respectively
\begin{equation}
\left| \int_{0}^T \int_{\R^3} A B u \bar v dx dt \right| \lessa
T^\delta \|
  B\|_{U^2_W L^2} \|
  A\|_{U^2_W H^1} \| u\|_{\tU^2_A H^2} \| v\|_{V^2_A L^2}
\label{h2b}\end{equation}

To prove \eqref{h2a} we use a triple Littlewood-Paley decomposition to
write
\[
\begin{split}
\left| \int_{0}^T \int_{\R^3} B \nabla u \bar v dx dt \right|
 \les & \sum_{\mu \les \l}
I^T_{\l,\l,\mu}(\nabla u,v,B) 
+\sum_{\mu \ll \l} I^T_{\mu,\l,\l}(\nabla u,v,B)
\\ & + \sum_{\mu \ll  \l}  I^T_{\l,\mu,\l}(\nabla u,v,B)
\end{split}
\]
Then for each term we use the corresponding bounds \eqref{llmh2h2},
and \eqref{lmlh2h2}   with a $\ln{\l}$ correction coming
from the use of Proposition \ref{principle}.   The summation
with respect to $\mu$ and $\l$ is straightforward.

For \eqref{h2b}, by \eqref{wcloc} the norms of $A$ and $B$ are $l^2$ summable 
with respect to unit spatial cubes. Hence without any restriction in
generality we can assume that both $A$ and $B$ are supported in a 
unit cube $Q$. For $A$ we use a Strichartz estimate, for $B$  the 
energy and  for $u$ a pointwise bound. Finally, for $v$ we interpolate
the local energy estimates with the local Strichartz estimates
to obtain
\[
\| S_\lambda v\|_{L^2  L^\frac{12}5} \lessa \lambda^{-\frac14+\e}
\|v\|_{V^2_A L^2}
\]
which leads to 
\[
\| v\|_{L^2  L^\frac{12}5} \lessa
\|v\|_{V^2_A L^2}
\]
Then we can estimate
\[
\begin{split}
  \left| \int_{0}^T \int_{Q} A B u \bar v dx dt \right| & \lesssim
  T^\frac14 \| B\|_{L^\infty L^2} \| A\|_{L^4 L^{12}} \|
  u\|_{L^\infty} \| v\|_{L^2 L^\frac{12}5} \\ & \lessa T^\frac14 \|
  B\|_{U^2_W L^2} \| A\|_{U^2_W H^1} \| u\|_{U^2_A H^2} \| v\|_{V^2_A
    L^2}
\end{split}
\]
\end{proof}

We now return to the equation \eqref{eqd}.  By 
the definition of the $\tU^2_AL^2$ norm and \eqref{atbd} we have 
\begin{equation}
\| 2 ( A_t \nabla -i A A_t ) ( 1-\Delta_A)^{-1} v\|_{DU^2_A L^2}
\lessa T^\delta \|v\|_{U^2_A L^2}
\label{pert}\end{equation}
By \eqref{u2solve} it follows that  we can solve 
\eqref{eqd} perturbatively in $U^2_A L^2$ on short time intervals.
This gives a solution $u = ( 1-\Delta_A)^{-1} v \in
\tU^2_A H^2$ for \eqref{schrainhom}. 
Furthermore, we obtain  the bound
\[
\| v(t) - S(t,0) v(0)\|_{L^2} \lessa T^\delta \|v(0)\|_{L^2}
\]
where $S(t,s)$ is the evolution associated to \eqref{schra}.
In particular this shows that $v \in C([0,1]; L^2)$.  An elliptic
argument allows us to return to $u$ and conclude that 
$u \in C([0,1]; H^2)$. This concludes the proof of the $H^2$
well-posedness.

Finally, we show that $U^2_A H^2 = \tU^2_A H^2$ and the other two
similar identities. Via the operator $I - \Delta_A$ these two spaces 
can be identified with the $ U^2 L^2$ spaces associated to 
the equations \eqref{schra}, respectively \eqref{eqd}.
But by Lemma~\ref{flowdiff}, these are equivalent due to \eqref{pert}.

\end{proof}

Next we consider the well-posedness in $H^{-2}$, which is essentially
dual to the $H^2$ well-posedness.

\begin{p1}
  Let $A \in U^2_W H^1$ with $ \nabla \cdot A = 0$. Then the equation
  \eqref{schrainhom} is well-posed in $H^{-2}$. In addition we have
\[
U^2_A H^{-2} = \tU^2_A H^{-2}, \qquad V^2_A H^{-2} = \tV^2_A H^{-2},
\qquad  DU^2_A H^{-2} = \tU^2_A H^{-2}
\]
with equivalent norms.
\end{p1}
\begin{proof}
By  Lemma~\ref{1-d}  we can write  the 
initial data $u_0$ as 
\[
u_0 = (1- \Delta_A) v_{0}, \qquad v_0 \in L^2
\]
Then we seek the solution $u$ for \eqref{schra} of the form $u =
(1-\Delta_A)v$.  The equation for $v$ is
\begin{equation}
(i \partial_t - \Delta_A) v = 2 (1-\Delta_A)^{-1}( A_t \nabla -i AA_t) v 
\label{perta}\end{equation}
To solve it we need the following counterpart to \eqref{atbd}:

\begin{l1}
 The operator $A_t \nabla - i A A_t$ satisfies the space-time bound 
\begin{equation}
\|  (A_t \nabla  - 2i A A_t) u \|_{DU^2_A H^{-2}} \lessa T^\delta \| u\|_{\tU^2_A L^2}
\label{atbd1}\end{equation}

\end{l1}
\begin{proof}
 By duality, \eqref{atbd1}  follows from the bounds
\begin{equation}
\left| \int_{0}^T \int_{\R^3} B \nabla u \bar v  dx dt \right| \lessa T^\delta
 \| u\|_{U^2_A L^2}  \| v\|_{V^2_A H^2}\|  B\|_{U^2_W L^2}
\label{h2a2} \end{equation}
respectively
\begin{equation}
\left| \int_{0}^T \int_{\R^3} A B u \bar v dx dt \right| \lessa
T^\delta \|
  B\|_{U^2_W L^2} \|
  A\|_{U^2_W H^1} \| u\|_{U^2_A L^2} \| v\|_{V^2_A H^2}
\label{h2b2}\end{equation}
These are almost identical to \eqref{h2a} and \eqref{h2b}, and their
proofs are essentially the same.
\end{proof}

The bound \eqref{atbd1} allows us to solve the equation \eqref{perta}
perturbatively in $U^2_A L^2$ and obtain a solution $v \in
C([0,1];L^2)$. This implies the $H^{-2}$ solvability for
\eqref{schra}. The second part of the proposition follows again from
Lemma~\ref{flowdiff}.
\end{proof}

Having the well-posedness result in $H^2$ and $H^{-2}$ allows 
us to prove well-posedness in a range of intermediate spaces.
Given a positive sequence $\{m(\lambda) \}_{\lambda = 2^j}$ 
satisfying 
\begin{equation}
0 < c < \frac{m(2\lambda)}{m(\lambda)} < C
\label{mfirst} \end{equation}
we define the Sobolev type space $H(m)$ with norm
\[
\| u\|_{H(m)}^2 = \sum_\lambda m^2(\lambda)
\|S_\lambda u\|_{L^2}^2 
\]
The standard Sobolev spaces $H^\alpha$ are obtained by taking
$m(\l)=\l^\alpha$. 

We consider the solvability for \eqref{schra} in $H(m)$ under a
stronger condition for $m$, namely
\begin{equation}
\frac14  \leq \frac{m(2\lambda)}{m(\lambda)} \leq 4
\label{alphaseq}\end{equation}
This guarantees that $H(m)$ is an intermediate space
between $H^{-2}$ and $H^2$.  

We need to describe $H(m)$ in terms of $H^{-2}$ and $H^2$. To measure
functions which are localized at some frequency $\lambda$ we 
can use the norm 
\[
\| u\|_{H_\lambda}^2 =  \lambda^{-4} \|u\|_{H^2}^2 +  \lambda^{4}
 \|u\|_{H^{-2}}^2
\]
Ideally we would like to represent $H(m)$ as an almost orthogonal
superposition of the $H_{\lambda}$ spaces with the weights
$m(\lambda)$. However, this does not work so well if $H(m)$ is
``close'' to either $H^{-2}$ or $H^2$.  Instead we need to select a
subset of dyadic frequencies which achieves the desired result.
We denote
\[
m_\infty = \lim_{\lambda \to \infty} \lambda^{2} m(\lambda)
\]
On $2^\N \cup \{ \infty \}$ we introduce the relation ``$\prec$'' by
\[
\lambda \prec \mu \Leftrightarrow  2 m(\mu) \geq  m(\lambda)
 (\lambda^2\mu^{-2}+ \mu^2 \lambda^{-2}), \qquad \mu, \lambda < \infty 
\]
respectively 
\[
\infty \prec \mu  \Leftrightarrow     2m(\mu) \geq  \mu^{-2} m_\infty,
\qquad \mu < \infty 
\]

\begin{d1} We say that a subset $\Lambda(m) \subset 2^\N \cup \{
  \infty \}$ is $m$-representative if (i) for each $\mu \in 2^\N$
  there exists $\lambda \in \Lambda(m)$ so that $ \lambda \prec \mu$
  and (ii) for each $\mu \in \Lambda(m) $ there is at most one $\lambda
  \in \Lambda(m) \setminus \{ \mu \}$ such that $\l  \prec \mu$.
\label{lofm}\end{d1}

\begin{l1}
  If $m$ satisfies \eqref{alphaseq} then an $m$-representative set
  $\Lambda(m)$ exists. In addition, for each $\mu \in 2^\N$ and $K
  \in \N$ we have
\begin{equation}
| \{ \l \in \Lambda(m);    2^K m(\mu) \geq m(\l)
 (\lambda^2\mu^{-2}+\mu^2 \lambda^{-2}) \} | \leq 4(K+4)
\label{lllb}\end{equation}
\end{l1}

\begin{proof}
  For each $\l \in 2^\N \cup \{ \infty \}$ we denote
\[
I_\l = \{ \mu \in 2^\N \cup \{ \infty \}, \ \l  \prec
\mu \}
\]
Due to \eqref{alphaseq} it is easy to see that $I_\l$ is an interval,
\[
I_\l = [ \l^-,\l^+], \qquad  \l^- \leq \l \leq \l^+
\]
and the endpoints $\l^-$ and $\l^+$ are nondecreasing functions of 
$\l$.  

We construct the set $\Lambda(m)$ as an increasing sequence
$\{\l_j\}$ in an iterative manner. $\l_0$ is chosen maximal 
so that $\l_0 \prec 1$. Iteratively, $\l_{j+1}$ is chosen
maximal so that $\l_{j+1} \prec 2 \l_j^+$. Either this process 
continues for an infinite number of steps, or it stops at some step
$k$ with $\l_k = \infty$. The former occurs if $m_\infty = \infty$
and the latter if $m_\infty < \infty$.

The property (i) in Definition~\ref{lofm} is satisfied by
construction. For (ii) we observe that $\l_{j+1} \geq 2 \l_j^+$
therefore $\l_j \not \prec \l_{j+1}$. On the other hand by
construction we have $\l_{j+2} \not\prec \l_j$.

For \eqref{lllb} suppose $\mu \leq \lambda_j$. Then
\[
m(\mu) \geq \frac{1}{4} m(\l_j) \mu^2 \l_j^{-2} \geq \frac{1}{8} m(\l_{j+2})
\mu^2 \l_{j+2}^{-2} \geq \cdots \geq 2^{-K-2} m(\l_{j+2K})
\mu^2 \l_{j+2K}^{-2} 
\]
A similar bound holds if we descend from $\mu$, and the conclusion
follows.
\end{proof}

Since we allow $\infty \in \Lambda(m)$ we need the equivalent of
$H_\l$ in that case, which is defined by $H_{\infty}=H^{-2}$. We note
that at  the other extreme we have  $H_1= H^2$.

\begin{l1} 
Let $m$ satisfy \eqref{alphaseq}, and $\Lambda(m)$ be an
$m$-representative subset of $2^N \cup \{\infty\}$.
Then 
\[
\| u\|_{H(m)}^2 \approx \inf \{ \sum_{\lambda \in \Lambda(m)}
m(\lambda)^2 \| u_\lambda \|_{H_\lambda}^2, \ u = \sum_{\lambda \in
  \Lambda(m)} u_\lambda \}
\]
\end{l1}

\begin{proof}
By Definition \ref{lofm} we have a finite  covering of $2^{\N}$ with intervals
\[
2^\N \subset \bigcup_{\lambda \in \Lambda(m)} I_\lambda
\]
We consider an associated partition of unity in the Fourier space,
\[
1 = \sum_{\lambda \in \Lambda(m)} \chi_\lambda(\xi)
\]
For $\mu \in I_\lambda$ we have $m(\mu) \approx m(\lambda)(\mu^2
\lambda^{-2} + \lambda^2 \mu^{-2})$ therefore we obtain
\[
\| \chi_\lambda(D_x) u\|_{H(m)} \approx m(\lambda)
\| \chi_\lambda(D_x) u\|_{H_\lambda} 
\]
and the ``$ \gtrsim $'' inequality follows. For the reverse we use 
\eqref{lllb}, which shows that the series $\sum u_\lambda$ is almost
orthogonal in $H(m)$,
\[
\langle u_{\lambda_i} , u_{\lambda_j} \rangle_{H(m)} \les 2^{-|i-j|}
 m(\lambda_i) m(\l_j)
\| u_{\lambda_i} \|_{H_{\lambda_i}} \|u_{\l_j}\|_{H_{\l_j}}
\]
\end{proof}

Finally we consider the well-posedness of \eqref{schra} in $H(m)$.

\begin{p1}
  a) Assume that the sequence $m$ satisfies \eqref{alphaseq}. Then the
  equation \eqref{schra} is well-posed in $H(m)$.
 
b)  Furthermore, for each $u \in U^2_A H(m)$ there is a representation
\[
u = \sum_{\lambda \in \Lambda(m)} u_\l
\]
 with 
\begin{equation}
 \sum_{\lambda \in \Lambda(m) } m^2(\lambda) 
\left( \lambda^{-4} \|u_\l \|_{U^2_A H^2}^2 +  \lambda^{4} \|u_\l \|_{U^2_A
  H^{-2}}^2\right) \lessa \| u\|_{U^2_A H(m)}^2
\label{u2hmrep}\end{equation}

c) The following duality relation holds:
\[
(DU^2_A H(m))^* = V^2_A H(m^{-1})
\]
\label{phm}\end{p1}

\begin{proof}
a) We consider a dyadic decomposition of the initial data 
\[
u_0 = \sum_ {\lambda \in \Lambda(m)} \chi_\lambda(D_x)  u_0
\]
and denote by $u_\lambda$ the solutions to \eqref{schra} with 
initial data $S_\lambda u_0$. Then
\[
u = \sum_ {\lambda \in \Lambda(m)}  u_\l
\]
We can measure $u_\l$ in both $H^2$ and $H^{-2}$,
\[
\l^{-2} \| u_\l \|_{C([0,1],H^2)} + \l^{2} \| u_\l
\|_{C([0,1],H^{-2})} \lessa \|  S_\lambda u_0\|_{L^2}
\]
After summation this gives
\[
 \sum_{\lambda\in \Lambda(m)} m^2(\lambda)
\left(\l^{-4} \| u_\l \|_{C([0,1],H^2)} + \l^{4} \| u_\l
\|_{C([0,1],H^{-2})} \right) \lessa \|u_0\|_{H(m)}^2
\]
which, by \eqref{alphaseq}, implies that $u \in C([0,1],H(m))$ and 
\[
\| u\|_{C([0,1],H(m))} \lessa \|u_0\|_{H(m)}
\]

b) It suffices to consider the case when $u$ is an $U^2_A H(m)$
atom. Then we consider a decomposition as in part (a) for
each of the steps, and the conclusion follows.

c) This is a direct consequence of \eqref{dual1}.
\end{proof}

Finally, we can interpolate the properties from $U^2_A H^2$ and $U^2_A
H^{-2}$ to obtain properties for $U^2H(m)$. Indeed, we have the following

\begin{p1} \label{prophm} Let $A \in U^2_W H^1$ with $ \nabla \cdot A
  = 0$. Assume that $m$, $m_1$, $m_2$  satisfy \eqref{alphaseq}.  Then
the following estimates hold:

(i) For $p,q$ as in \eqref{seu2} we have the Strichartz estimate
\begin{equation} \label{hastr} \| S_\lambda u \|_{L^p(0,T; L^q)} \lessa
  T^{\frac{\delta}p} m(\l)^{-1} \lambda^{\frac1p} \| u\|_{U^2_A H(m)}
\end{equation}

(ii) Local energy estimates.  For any spatial cube $Q$ of size  $1$ we have:
\begin{equation}
  \| S_\lambda u\|_{L^2([0,T]\times Q)} \lessa  T^\delta m(\lambda)^{-1}
  \lambda^{- \frac12 +\epsilon}  \| u_0\|_{\tU^2_A H(m)}
\label{lehm}\end{equation} 

(iii)  Local Strichartz estimates. For any spatial cube $Q$ of size  $1$ we have:
\begin{equation}
  \| S_\lambda u\|_{L^2(0,1; L^6(B))}  \lessa  T^\delta m(\lambda)^{-1} 
\lambda^{\epsilon} \| u_0\|_{U^2_A H(m)}
\label{seulochm}\end{equation} 

(iv) Trilinear estimates. For $\mu \leq \lambda$ we have
\begin{equation}
\!\! \!\! |I^T_{\l,\l,\mu}(u,v,B)| \! \lessa \! \frac{T^{\delta} \l^\epsilon \min{(1,\mu \l^{-\frac12})}}{ m_1(\l) m_2(\l)} \|u\|_{U^2_A H(m_1)}  \|v\|_{U^2_A H(m_2)}
\|B\|_{U^2_W L^2}  
\label{llmha}\end{equation}
while if $\mu \ll \lambda$ then
\begin{equation}
|I^T_{\l,\mu,\l}(u,v,B)| \lessa \frac{T^\delta 
\lambda^{-\frac12+\epsilon} \mu^\frac12}{m_1(\l) m_2(\mu)}  \|u\|_{U^2_A H(m_1)}  \|v\|_{U^2_A H(m_2)}
\|B\|_{U^2_W L^2}  
\label{lmlha}\end{equation}
\end{p1}

Due to the representation in Proposition~\ref{phm}(b) this result is a
straightforward consequence of the similar results in $H^2$ and $H^{-2}$.

\section{ Proof of Theorem~\ref{tmain}}
\label{stmain}

We first establish an a-priori estimate for regular ($H^2$) solutions of 
\eqref{msc}. For this we need to consider the two nonlinear
expressions on the right hand side of \eqref{msc}. We begin with the
Schr\"odinger nonlinearity:

\begin{l1} \label{schnonlin}
For $\beta > \frac12$ and $m$ satisfying \eqref{alphaseq}  we have 
\begin{equation}
\| \phi v\|_{DU^2_A(0,T; H(m))} \lessa T^\delta \| v\|_{U^2_A H(m)}
\| u\|_{U^2_A H^\beta}^2
\label{phiu}\end{equation}

\end{l1}
\begin{proof}
  By duality the above bound is equivalent to the quadrilinear
  estimate
\[
\left|\!\inT \Delta^{\!-1} \!(u_1 \bar u_2) u_3 \bar u_4 dx dt \right| \!
\lessa \! T^\delta \| u_1\|_{U^2_A H^\beta}\| u_2\|_{U^2_A H^\beta}
\| u_3\|_{U^2_A H(m)}\| u_4\|_{V^2_A H(\frac{1}m)}
\]
After a simultaneous Littlewood-Paley decomposition of the three factors
$\Delta^{-1} (u_1 \bar u_2)$, $u_3$ and $u_4$ we need to consider
the following three sums:
\[
S_{lhh} =  \sum_{\mu \leq \lambda} \inT \Delta^{-1} S_\mu  (u_1 \bar u_2)
S_\lambda u_3 S_\lambda \bar u_4 dx dt 
\]
\[
S_{hlh} = \sum_{\mu \ll \lambda} \inT \Delta^{-1} S_\lambda (u_1 \bar
u_2) S_\mu u_3 S_\lambda \bar u_4 dx dt
\]
\[ 
S_{hhl} = \sum_{\mu \ll \lambda} \inT \Delta^{-1} S_\lambda (u_1 \bar
u_2) S_\l u_3 S_\mu \bar u_4 dx dt
\]

The first sum can be estimated using Strichartz and energy estimates
as follows:
\[
\begin{split}
|S_{lhh}| \les & \sum_{\mu} \mu^{-2} \| S_\mu  (u_1 \bar u_2) \|_{L^1 L^\infty} 
\Big\| \sum_{\l \geq \mu} S_\lambda u_3  S_\lambda \bar u_4  
\Big\|_{L^{\infty} L^1}
\\ 
\les & \sum_{\mu} \mu^{-1} \| S_\mu  (u_1 \bar u_2) \|_{L^1 L^3}
\sup_{t \in [0,T]} \sum_{\l \geq \mu} \|S_\l u_3(t)\|_{L^2} \|S_\l u_4(t)\|_{L^2}
\\
\les & \| u_1 \|_{L^2 L^6} \| u_2 \|_{L^2 L^6} 
\sup_{t \in [0,T]}  \|u_3(t)\|_{H(m)} \|u_4(t)\|_{H(m^{-1})}
\\
\les_A & T^\delta \| u_1\|_{U^2_A H^\beta}\| u_2\|_{U^2_A H^\beta}
\| u_3\|_{U^2_A H(m)}\| u_4\|_{V^2_A H(m^{-1})}
\end{split}
\]

The second and third sums are similar. Using the Strichartz 
estimates we obtain
\[
\begin{split}
  |S_{hlh}| \les &  \sum_{\mu \ll \lambda} \l^{-2} \| S_\l (u_1
  \bar u_2) \|_{L^2 L^2} \| S_\mu u_3 \|_{L^2 L^\infty} \| S_\lambda
  \bar u_4 \|_{L^{\infty} L^2}
  \\
  \les_A & T^{\delta}  \sum_{\mu \ll \lambda}
   \frac{\mu m(\l)}{\l^\frac32 m(\mu)} \| S_\l (u_1 \bar
  u_2) \|_{L^2 L^\frac32} \| u_3 \|_{U^2_A H(m)} \| u_4 \|_{V^2_A H(m^{-1})}
\end{split}
\]
At least one of $u_1$ or $u_2$, say $u_1$, must have  frequency
at least $\lambda$. Then we continue  with
\[
\begin{split}
|S_{hlh}| \les_A & T^{{\delta}}  \sum_{\mu \ll \lambda} \frac{\mu
  m(\l)}{\l^\frac32 m(\mu)} \| u_1 \|_{L^\infty L^2} \| \bar u_2
\|_{L^2 L^6} \| u_3 \|_{U^2_A H(m)} \| u_4 \|_{V^2_A H(m^{-1})}
\\
\les_A & T^{\delta} \sum_{\mu \ll \lambda}
\frac{\mu^2m(\l)}{\l^2 m(\mu)}\frac{ \l^{\frac12 - \beta}} {\mu} \| u_1
\|_{U^2_A H^\beta} 
\| \bar u_2 \|_{U^2_A H^\beta} \| u_3 \|_{U^2_A H(m)} \|  u_4 \|_{V^2_A H(m^{-1})}
\end{split}
\]
By \eqref{alphaseq}  the first fraction is less than one therefore the
summation with respect to $\lambda$ and $\mu$ is straightforward. 
 \end{proof}

 Next we consider the wave nonlinearity. If $m$ is as in
 \eqref{mfirst} then the linear wave equation is well-posed in $H(m)$,
 and can easily define the corresponding spaces $U^2_W H(m)$, $V^2_W
 H(m)$ respectively $DU^2_W H(\lambda^{-1}m)$.  The next result
asserts that in effect the contribution of the wave nonlinearity is
one half of a derivative better than the solution to the Schr\"odinger
equation. Here we impose an additional condition on $m$, namely
\begin{equation}
\frac{m(\lambda)}{m(\mu)} \geq \frac{\lambda^\beta}{\mu^\beta}, \qquad
\lambda > \mu
\label{mbhigh}\end{equation}
which guarantees that the $H(m)$ norm is at least as strong as the
$H^\beta$ norm.

\begin{l1}\label{wnonlin}
  a) Let $\beta > \frac12$ and $m$ satisfying \eqref{alphaseq} and
  \eqref{mbhigh}. Then
\begin{equation}
\begin{split}
\| P(\bar u \nabla_A u)\|_{DU^2_W(0,T; H(\lambda^{-\frac12} m))} 
\! \lessa &\  T^\delta( \| u\|_{U^2_A H(m)}  \| u\|_{U^2_A H^\beta}  
\\ &\
 +
\| u\|_{U^2_A H^\beta}^2 \|A\|_{U^2_W H(m)})
\end{split}
\label{udau}\end{equation}

b) For $\beta > \frac34$ we have
\begin{equation}
\|P(\bar u \nabla_A v)\|_{DU^2_W(0,T; H^{-\frac12})} 
\lessa T^\delta \| u\|_{U^2_A L^2}  \| v \|_{U^2_A H^\beta}  
\label{udaucd}\end{equation}

\end{l1}
\begin{proof}
a) By duality we have two estimates to prove. The first is
\begin{equation}
\left|  \inT \bar u \nabla u B dx dt \right| 
\lessa T^\delta \| u\|_{U^2_A  H(m)}  
\| u\|_{U^2_A H^\beta}  \| B\|_{V^2_W H(\lambda^\frac12/ m)}
\label{nonlinw3} \end{equation}
with divergence free $B$. The second is
\begin{equation}
\begin{split}
\left|  \inT  \bar u  u  A B dx dt \right| & \lessa  T^\delta \| u\|_{U^2_A
  H(m)}  \| u\|_{U^2_A H^\beta}  \| A\|_{V^2_W H^1}
 \| B\|_{V^2_W H(\lambda^\frac12/ m)}
\\
& + T^\delta \| u\|_{U^2_A
  H^\beta}  \| u\|_{U^2_A H^\beta}  \| A\|_{V^2_W H(m)}
 \| B\|_{V^2_W H(\lambda^\frac12/ m)}
\end{split}
\label{nonlinw4} 
\end{equation}

Consider \eqref{nonlinw3}. Since $B$ is divergence free, it follows 
that the gradient can be placed either on $u$ or on $\bar u$. 
Hence using a simultaneous trilinear Littlewood-Paley
decomposition of the three factors we reduce the problem to
estimating the following two terms:
\[
S_{hhl} = \sum_\l \sum_{\mu \les \l}| I^T_{\l,\l,\mu}( u,\nabla
u,B)|, \qquad S_{hlh} = 
\sum_\l \sum_{\mu \ll \l} | I^T_{\l,\mu,\l}(u,\nabla
u,B)|
\]
We use \eqref{llmha} and Lemma~\ref{principle} to estimate the first term:
\[
S_{hhl} \lessa T^{\delta} \sum_\l \sum_{\mu \les \l} \frac{
  \l^{\epsilon} m(\mu)}
 {\l^{\beta} \mu^{\frac12} m(\l)} 
\| u\|_{U^2_A  H(m)}  
\| u\|_{U^2_A H^\beta}  \| B\|_{V^2_W H(\lambda^\frac12/ m)}
\]
The bound \eqref{mbhigh} insures that the summation is straightforward
if $\e$ is chosen sufficiently small.
For the second term we use \eqref{lmlha} and  Lemma~\ref{principle}:
\[
S_{hlh} \lessa T^{\delta} \sum_\l \sum_{\mu \les \l} \frac{
  \mu^{\frac32-\beta}}{\lambda^{1-\e}} \| u\|_{U^2_A H(m)} \|
u\|_{U^2_A H^\beta} \| B\|_{V^2_W H(\lambda^\frac12/ m)}
\]
which is again summable if $\e$ is sufficiently small.

Next we turn our attention to \eqref{nonlinw4}. 
Using a Littlewood-Paley decomposition for all terms
 it suffices to consider factors of type
\[
I_{\lambda_1,\lambda_2,\lambda_3,\lambda_4}^T =\left| \int_0^T
  \int_{\R^3}S_{\lambda_1}  \bar  u \, S_{\lambda_2} u\,
  S_{\lambda_3} A \, S_{\lambda_4} B dx dt \right|
\]
and prove that for some $\delta > 0$ they satisfy the bound
\begin{equation}
I_{\lambda_1,\lambda_2,\lambda_3,\lambda_4}^T \lessa \lambda_1^{-\delta}
 \lambda_2^{-\delta} \lambda_3^{-\delta} \lambda_4^{-\delta}
 \cdot RHS(\eqref{nonlinw4})
\label{es4w}\end{equation}
We begin with a weaker bound which follows directly from 
Strichartz estimates, namely
\[
I_{\lambda_1,\lambda_2,\lambda_3,\lambda_4}^T \lessa T^\frac13
\lambda_1^{-\frac14} \lambda_2^{-\frac14} \lambda_3^{-\frac23}
\lambda_4^N \|u\|_{U^2_A H^\beta}^2 \|A\|_{U^2_W H^1} \| B\|_{L^\infty
  L^2}
\]
Arguing as in Lemma~\ref{principle}, this allows us to replace \eqref{es4w}
with 
\begin{equation}
I_{\lambda_1,\lambda_2,\lambda_3,\lambda_4}^T \lessa \lambda_1^{-\delta}
 \lambda_2^{-\delta} \lambda_3^{-\delta} \lambda_4^{-\delta}
 \cdot {\text{ modified } RHS}(\eqref{nonlinw4})
\label{es4wa}\end{equation}
where we have replaced $V^2_W H(\lambda^\frac12/m)$ space with the
similar $U^2$ space in the right hand side of \eqref{nonlinw4}.

Due to \eqref{wcloc} both wave
factors are $l^2$ summable with respect to unit spatial cubes
therefore it is enough to estimate the above integral on a unit cube
$Q$.  Also we must have $\lambda_4 \leq \max
\{\lambda_1,\lambda_2,\lambda_3\}$. Hence we consider two cases.  

If $\max\{\lambda_1,\lambda_2,\lambda_3\}=\l_1$ (or $\lambda_2$) then we
estimate
\[
\begin{split}
& I_{\lambda_1,\lambda_2,\lambda_3,\lambda_4}^T 
\leq T^\frac16 \| S_{\lambda_1} u\|_{L^4 L^3} \|S_{\lambda_2} u\|_{L^2 L^\infty}
 \|S_{\lambda_3} A\|_{L^3 L^6}    \|S_{\lambda_4} B\|_{L^\infty L^2}  
\\ 
\leq & \  T^\frac16 \lambda_1^{-\frac12+\e} 
 \lambda_2^{\frac12+\e-\beta} \lambda_3^{-\frac16}
 \frac{\lambda_1^\frac12 m(\lambda_4)}{\lambda_4^\frac12 m(\lambda_1)}
      \| u\|_{U^2_A
  H(m)}  \| u\|_{U^2_A H^\beta}  \| A\|_{V^2_W H^1}
 \| B\|_{U^2_W H(\frac{\sqrt\lambda}{ m})}
\end{split}
\]
By \eqref{mbhigh} the fraction above is less than one therefore for
small enough $\e$ the bound \eqref{es4wa} follows.  

The second case is when $\max\{\lambda_1,\lambda_2,\lambda_3\}=\l_3$. Then we estimate
\[
\begin{split}
& I_{\lambda_1,\lambda_2,\lambda_3,\lambda_4}^T 
\leq T^\frac1{12} \| S_{\lambda_1} u\|_{L^4 L^6} \|S_{\lambda_2} u\|_{L^2 L^\infty}
 \|S_{\lambda_3} A\|_{L^6 L^3}    \|S_{\lambda_4} B\|_{L^\infty L^2}  
\\ 
\lessa &   T^\frac1{12} \lambda_1^{\frac13+\e-\beta}  
 \lambda_2^{\frac12+\e-\beta} \lambda_3^{-\frac16}
 \frac{\lambda_3^\frac12 m(\lambda_4)}{\lambda_4^\frac12 m(\lambda_3)}
      \| u\|_{U^2_A
  H^\beta}  \| u\|_{U^2_A H^\beta}  \| A\|_{U^2_W H(m)}
 \| B\|_{U^2_W H(\frac{\sqrt\lambda}{ m})}
\end{split}
\]
and \eqref{es4wa} again follows.

b) By duality we have two estimates to prove. The first one is trilinear,
\begin{equation}
\left|  \int_0^T \int_{\R^3} \bar u \nabla v B dx dt \right| 
\lessa T^\delta \| u \|_{U^2_A
  L^2}  \| v\|_{U^2_A H^\beta}  \| B\|_{V^2_W H^{\frac12}}
\label{nonlinw5} \end{equation}
with a divergence free $B$. The second one is quadrilinear,
\begin{equation}
\left|  \inT  \bar u  u  A B dx dt \right|  \lessa  T^\delta \| u\|_{U^2_A
 H^\beta }  \| v \|_{U^2_A L^2}  \| A\|_{U^2_W H^1}
 \| B\|_{V^2_W H^\frac12}
\label{nonlinw6} 
\end{equation}

Since $B$ is divergence free, in \eqref{nonlinw5} we can place the
gradient on either $u$ or $v$. The argument is similar to the one in
part (a), using \eqref{llmha} and \eqref{lmlha}, as well as
Lemma~\ref{principle} in order to substitute the $V^2$ norm by the
$U^2$ norm. The new restriction $\beta > \frac34$ arises in the case
when the low frequency is on $B$. Indeed, if $\mu \ll \lambda$ then by
\eqref{llmha} we have
\[
|I^T_{\l,\l,\mu} (u,\nabla v,B)|\! \lessa \! T^\delta
\lambda^{1-\beta+\epsilon} \mu^{-\frac12}
 \min\{ 1, \mu \lambda^{-\frac12}\}
 \| u \|_{U^2_A
  L^2}  \| v\|_{U^2_A H^\beta}  \| B\|_{V^2_W H^{\frac12}}
\]
The worst case is when $\mu =\lambda^\frac12$, when the 
coefficient above is $T^\delta \lambda^{\frac34-\beta+\e}$.

The estimate \eqref{nonlinw6} is also proved as in part (a). Indeed,
by Corollary~\ref{chhk} we can substitute the $V^2$ space by $U^2$ at
the expense of losing $\e$ derivatives. Because of the finite speed of
propagation for the wave equation, see \eqref{wcloc}, we can reduce
the problem to the case when $A,B$ are supported in a unit cube $Q$.
There we can use the Strichartz estimates for the wave equation,
respectively the local Strichartz estimates for the Schr\"odinger
equation.

\end{proof}

The next step in the proof of the theorem is to establish an
a-priori $H^1$ estimate for $H^2$ solutions. This is obtained 
in terms of the conserved quantities in our problem, namely $E$ and $Q$.

\begin{p1} 
Let $(u,A)$ be an $H^2$ solution for \eqref{msc} in some 
time interval $[0,T_0]$ with $T_0 \leq 1$. Then 
\[
\| u\|_{U^2_A (0,T_0;H^1)} + \| A \|_{U^2_W(0,T_0;H^1)} 
\leq c(E,Q)
\]
\label{ph1}
\end{p1}

\begin{proof}
We use a bootstrap argument. Since $u,A \in L^\infty H^2$,
we can easily estimate the wave nonlinearity and obtain 
$\Box A \in L^\infty H^1$. This implies that $A \in U^2_WH^1$,
and, in addition, that the function
\[
T \to \|A\|_{U^2_W(0,T;H^1)} 
\]
is continuous and satisfies
\[
\lim_{T \to 0} \|A\|_{U^2_W(0,T;H^1)} = \|A(0)\|_{H^1}+\|A_t(0)\|_{L^2}
\]
A similar argument applies in the case of the Schr\"odinger
equation.

By \eqref{u2solve}  we estimate in the  Schr\"odinger 
equation
\[
\| u \|_{U^2_A(0,T; H^1)} \lesssim
\|u(0)\|_{H^1} +   \| iu_t -\Delta_A u \|_{DU^2_A(0,T; H^1)}.
\]
Then by  \eqref{phiu} we obtain
\[
\| u \|_{U^2_A(0,T; H^1)} \leq C_A^1( \|u(0)\|_{H^1} +  
 T^\delta  \| u \|_{U^2_A(0,T; H^1)}^3) 
\]
Similarly we can use \eqref{udau} to obtain a bound for the wave 
equation
\[
\|A\|_{U^2_W(0,T;H^1)}  \leq  \|A(0)\|_{H^1}+\|A_t(0)\|_{L^2}+
T^\delta   C_A^2 \| u \|_{U^2_A(0,T; H^1)}^2
\]
We multiply the first equation by $c_A^1 = (C_A^1)^{-1}$ and add
to the second equation to obtain
\[
\begin{split}
\|A\|_{U^2_W(0,T;H^1)} + c_A^1 \| u \|_{U^2_A(0,T; H^1)} \leq 
 &\  
     \|u(0)\|_{H^1}+ \|A(0)\|_{H^1}+\|A_t(0)\|_{L^2}
\\ &\ + T^\delta C_A^3 (1 +\| u \|_{U^2_A(0,T; H^1)}^3) 
\end{split}
\]

We  make the bootstrap assumption
\[
c_A^1 \| u \|_{U^2_A(0,T; H^1)} + \|A\|_{U^2_W(0,T;H^1)} \leq 2 +
\|u(0)\|_{H^1}+ \|A(0)\|_{H^1}+\|A_t(0)\|_{L^2} 
\]
Then the previous bound implies that
\[
\begin{split}
c_A^1 \| u \|_{U^2_A(0,T; H^1)} + &\  \|A\|_{U^2_W(0,T;H^1)}
\leq 
     \|u(0)\|_{H^1}+ \|A(0)\|_{H^1}+\|A_t(0)\|_{L^2}
\\ &\ + T^\delta C(E,Q) 
(1+\| u \|_{U^2_A(0,T; H^1)}^3)
\end{split}
\]
This shows that for $T \leq T_0(E,Q)$ we have 
\begin{equation}
\| u \|_{U^2_A(0,T; H^1)} + \|A\|_{U^2_W(0,T;H^1)}
\leq  1+ \|u(0)\|_{H^1}+ \|A(0)\|_{H^1}+\|A_t(0)\|_{L^2}
\label{ste}\end{equation}
improving our bootstrap assumption. 

Hence a continuity argument shows that \eqref{ste}
holds without any bootstrap assumption. The conclusion of the 
proposition follows by summing up with respect to $T_0(E,Q)$
time intervals.

\end{proof}

The next step is to establish an a-priori $H^2$ bound
with constants which depend only on the $H^1$ size of the
data.

\begin{p1} 
Let $(u,A)$ be an $H^2$ solution for \eqref{msc} in some 
time interval $[0,T_0]$ with $T_0 \leq 1$. Then 
\[
\| u\|_{U^2_A H^2} + \| A \|_{U^2_WH^2} 
\leq c(E,Q) (\|u_0\|_{H^2} + \|A_0\|_{H^2} + \|A_1\|_{H^1})
\]
\end{p1}
 
\begin{proof}
The argument is similar to the one above.
\end{proof}

Given the local well-posedness result in $H^2$ proved in earlier work
\cite{MR2181763}, we can iterate the argument and conclude that the
$H^2$ solutions are global.

Finally, our last apriori estimate is in intermediate spaces:

\begin{p1} 
Let $m$ be a weight which satisfies \eqref{alphaseq} and \eqref{mbhigh}.
Let $(u,A)$ be an $H^2$ solution for \eqref{msc} in the
time interval $[0,1]$. Then 
\[
\| u\|_{U^2_A H(m)} + \| A \|_{U^2_WH(m)} 
\leq c(E,Q) (\|u_0\|_{H(m)} + \|A_0\|_{H(m)} + 
\|A_1\|_{H(\lambda^{-1}m)})
\]
\label{phalpha}
\end{p1}
\begin{proof}
The argument is again similar to the one above.
\end{proof}

In order to obtain $H^1$ solutions and to study the dependence
of the solutions on the initial data we need to obtain estimates 
for differences of solutions. Given a solution $(u,A)$ to \eqref{msc}
we consider the corresponding linearized problem
\begin{equation}
\left\{
\begin{array}{l} i v_t - \Delta_A v = 2iB \nabla u + 2 AB u + \phi v + 
\Delta^{-1}(\Re u \bar v) u \cr \cr
\Box B = P(\bar v \nabla_A u + \bar u \nabla_A v + B|u|^2)
\end{array} 
\right.
\label{linearized}
\end{equation}

Our main estimate for the linearized problem is 

\begin{p1}
Let $(u,A)$ be an $H^2$ solution for \eqref{msc} in the
time interval $[0,1]$. Then the linearized problem \eqref{linearized}
is well-posed in $L^2 \times H^\frac12 \times H^{-\frac12}$
uniformly with respect to $(u,A)$ in a bounded set in the energy 
space,
\begin{equation}
\| v\|_{U^2_A L^2} + \| B \|_{U^2_WH^\frac12} 
\leq c(E,Q) (\|v_0\|_{L^2} + \|B_0\|_{H^\frac12} + 
\|B_1\|_{H^{-\frac12}})
\end{equation}
\label{plinearized}\end{p1}

\begin{proof}
The conclusion follows iteratively in short time intervals provided 
that we obtain appropriate  estimates for the  terms on the right:
\[
\| 2iB \nabla u + 2 AB u + \phi v +
\Delta^{-1}(\Re u \bar v) u\|_{DU^2_A L^2} \! \lesssim T^\delta c(E,Q)
(\| v\|_{U^2_A L^2} +  \| B \|_{U^2_WH^\frac12})
\]
respectively
\[
\|P(\bar v \nabla_A u + \bar u \nabla_A v + iB|u|^2)\|_{DU^2_W
  H^{-\frac12}} 
 \lesssim T^\delta c(E,Q)
(\| v\|_{U^2_A L^2} +  \| B \|_{U^2_WH^\frac12})
\]
These in turn follow by duality from the trilinear and quadrilinear
bounds
\begin{equation}
\left| \inT  B \nabla u_1  \bar u_2 dx dt \right| \lessa  T^\delta
\| B \|_{U^2_WH^\frac12} \| u_1\|_{U^2_A H^1} \| u_2\|_{V^2_A L^2}
\label{lin1}\end{equation}
\begin{equation}
\left| \inT AB u_1 \bar u_2  dx dt \right| \lessa  T^\delta\| A \|_{U^2_WH^1}
\| B \|_{U^2_WH^\frac12} \| u_1\|_{U^2_A H^1} \| u_2\|_{V^2_A L^2}
\label{lin2}\end{equation}
\begin{equation}
\!\! \!\! \!\!
\left| \inT \Delta^{-1} (u_1 \bar u_2) u_3 \bar u_4  dx dt \right|\! \lessa \!T^\delta
  \| u_1\|_{U^2_A H^1} \| u_2\|_{U^2_A H^1} 
  \| u_3\|_{U^2_A L^2} \| u_4\|_{V^2_A L^2}
\label{lin3}\end{equation}
\begin{equation}
\!\!\! \! \!\!\left|  \inT \Delta^{-1} (u_1 \bar u_2) u_3 \bar u_4  dx dt
  \right|\! \lessa \! T^\delta
  \| u_1\|_{U^2_A H^1} \| u_2\|_{U^2_A L^2} 
  \| u_3\|_{U^2_A H^1} \| u_4\|_{V^2_A L^2}
\label{lin4}\end{equation}
\begin{equation}
\left|  \inT       B \nabla u_1  \bar u_2            dx dt \right| \lessa T^\delta
\| B \|_{V^2_W H^\frac12} \| u_1\|_{U^2_A H^1} \| u_2\|_{U^2_A L^2}
\label{lin5}\end{equation}
\begin{equation}
\left|  \inT     B \nabla u_1  \bar u_2            dx dt \right| \lessa T^\delta
\| B \|_{V^2_W H^\frac12} \| u_1\|_{U^2_A L^2} \| u_2\|_{U^2_A H^1}
\label{lin6}\end{equation}
\begin{equation}
\left| \inT   AB u_1 \bar u_2                 dx dt \right| \lessa T^\delta
\| A \|_{U^2_WH^1}
\| B \|_{V^2_WH^\frac12} \| u_1\|_{U^2_A H^1} \| u_2\|_{U^2_A L^2}
\label{lin7}\end{equation}
\begin{equation}
\left|  \inT   AB u_1 \bar u_2    dx dt \right| \lessa T^\delta
\| A \|_{U^2_WH^{\frac12}}
\| B \|_{V^2_WH^\frac12} \| u_1\|_{U^2_A H^1} \| u_2\|_{U^2_A H^1}
\label{lin8}\end{equation}

The quadrilinear mixed bounds \eqref{lin2}, \eqref{lin7}, \eqref{lin8}
follow trivially from the Strichartz estimates. For \eqref{lin2}
for instance we have 
\[
\begin{split}
\left| \inT AB u_1 \bar u_2  dx dt \right| \lesssim &\
T^{\frac{1}{12}} \| A \|_{L^3  L^8}
\| B \|_{L^4} \| u_1\|_{L^{3} L^8} \| u_2\|_{L^\infty L^2}
\\ \lessa &\ 
T^{\frac{1}{12}} \| A \|_{U^2_WH^1}
\| B \|_{U^2_WH^\frac12} \| u_1\|_{U^2_A H^1} \| u_2\|_{V^2_A L^2}
\end{split}
\]
We note that there is significant room for improvement in this
computation by localizing first to the unit spatial scale and then
using the local Strichartz estimates for the Schr\"odinger equation.

The quadrilinear Schr\"odinger bound \eqref{lin3} corresponds to
the particular choice $m(\l)=1$ and $\beta=1 > \frac12$ in \eqref{phiu}.
For \eqref{lin4} we can write
\[
\begin{split}
  \left| \inT \Delta^{-1} (u_1 \bar u_2) u_3 \bar u_4 dx dt \right|
  \lesssim &\ \| u_1 u_2\|_{L^2 L^\frac32} \| u_3 u_4\|_{L^2 L^\frac32}
  \\
  \lesssim &\ T^\frac13 \|u_1\|_{L^3 L^6} \|u_2\|_{L^\infty L^2}
  \|u_3\|_{L^3 L^6} \|u_4\|_{L^\infty L^2}
\\
\lessa &\ T^\frac13
  \| u_1\|_{U^2_A H^1} \| u_2\|_{U^2_A L^2} 
  \| u_3\|_{U^2_A H^1} \| u_4\|_{V^2_A L^2}
\end{split}
\]

Finally, the bounds \eqref{lin5} and \eqref{lin6} are identical since
$\text{div}\ B=0$ and correspond to \eqref{nonlinw5}. \eqref{lin1} is
essentially the same estimate.

% {\bf Note:} The short time apriori estimates work with $H^\beta$
% Schr\"odinger data for all $\beta > \frac12$. However, above it seems
% that we need a factor of $\mu^\frac12 \lambda^{\beta -1}$ in
% \eqref{llm}, which we get for $\beta = 1$ but, it seems, if $\mu =
% \lambda^\frac12$, we are limited to $\beta > 3/4$.
\end{proof}

{\bf Proof of Theorem~\ref{tmain}, conclusion:}

By Proposition~\ref{plinearized} we can obtain a weak Lipschitz
dependence result for $H^2$ solutions $(u_1,A_1)$ and $(u_2,A_2)$ to
\eqref{msc},
\begin{equation}
\begin{split}
\| u_1 -u_2\|_{L^\infty L^2} + \|A_1-A_2\|_{U^2_W H^\frac12}
\leq c(E_1,Q_1,E_2,Q_2)\\ ( \| (u_1 -u_2)(0)\|_{ L^2}
+ \|(A_1-A_2)(0)\|_{H^\frac12} + \|(A_1-A_2)_t(0)\|_{H^{-\frac12}})
\end{split}
\label{lowlip}\end{equation}

We use this in order to construct solutions to \eqref{msc} for $H^1$ 
initial data. Given $(u_0, A_0, A_1) \in H^1 \times H^1 \times L^2$
we consider a sequence of $H^2$ initial data
\[
(u_0^n,A_0^n,A_1^n) \to  (u_0, A_0, A_1) \text{ in }
 H^1 \times H^1 \times L^2
\]
The sequence $(u_0^n,A_0^n,A_1^n) $ is compact in  $H^1 \times H^1
\times L^2$, therefore we can bound them uniformly in a stronger norm,
\[
\| (u_0^n,A_0^n,A_1^n)\|_{H(m) \times H(m) \times H(\l^{-1}m)}
\leq M
\]
where $m(\l) \geq \l$ satisfies \eqref{alphaseq} and \eqref{mbhigh}
 and in addition
\[
\lim_{\lambda \to \infty} \lambda^{-1} m(\l) = \infty
\]
By Proposition~\ref{phalpha} we obtain a uniform bound
\[
\| u^n\|_{L^\infty H(m)} + \|A^n\|_{U^2_W H(m)}
\leq M
\]
On the other hand, \eqref{lowlip} shows that the solutions
$(u^n,A^n)$ have a limit in a weaker topology,
\[
(u^n,A^n) \to (u,A) \qquad \text{ in } L^\infty L^2 \times U^2_W H^\frac12
\]
Combining the two bounds above we obtain strong convergence in $H^1$,
\[
(u^n,A^n) \to (u,A) \qquad \text{ in } L^\infty H^1 \times U^2_W H^1
\]
In addition, $u$ will also satisfy the same Strichartz estimates as
$u^n$.  Passing to the limit in the equation \eqref{msc} we easily see
that $(u,A)$ is a solution.  Due to the weak Lipschitz dependence it
is also the unique uniform limit of strong solutions.  Due to the
Strichartz estimates we can bound the nonlinear term $\phi u$ in the
Schr\"odinger equation as in Lemma~\ref{phiu}. Then  it also follows that $u
\in U^2_A H^1$.

The weak Lipschitz dependence \eqref{lowlip} carries over to $H^1$
solutions, as well as the bounds in
Propositions~\ref{ph1},\ref{phalpha}. Then the same argument as above
gives the continuous dependence on the initial data.

\section{ Wave packets for Schr\"odinger operators with 
rough symbols}
\label{scale1}

An essential part of this article is devoted to understanding the
properties of the \eqref{schralambda} flow at frequency $\lambda$ on
$\lambda^{-1}$ time intervals.  As it turns out, for many estimates
the parameter $\lambda$ can be factored out by rescaling.
This is why in this section we consider a more general equation 
of the form
\begin{equation} \label{generic}
i u_t - \Delta u + a^w(t,x,D) u = 0, \qquad u(0) = u_0
\end{equation}
which we study on a unit time scale. Here $a$ is a  real
symbol which is roughly smooth on the unit scale.

For such a problem one seeks to obtain a wave packet
parametrix, i.e. to write solutions as almost orthogonal
superpositions of wave packets, where the wave packets 
are localized both in space and in frequency on the unit scale.
The simplest setup is to assume uniform bounds on $a$ of the  form
\[
| \partial_x^\alpha \partial_\xi^\beta a(t,x,\xi)| \leq c_{\alpha
  \beta}, \qquad |\alpha| + |\beta| \geq k
\]
An analysis of this type has been carried out in \cite{MR1887639},
\cite{MR2208883}, \cite{MR2094851}. If $k=2$ then one obtains a wave
packet parametrix where the packets travel along the Hamilton flow. If
$k=1$ the geometry simplifies, and the Hamilton flow stays close to
the flow for $a=0$; however, $a$ still affects a time modulation
factor arising in the solutions. Finally if $k=0$ then the
$a^w(t,x,D)$ term is purely perturbative.

For the operators arising in the present paper the above uniform
bounds on $a$ are too strong, and need to be replaced  by
integral bounds of the form
\[
\int_{0}^1 | \partial_x^\alpha \partial_\xi^\beta a(t,x^t,\xi^t)| \leq c_{\alpha
  \beta}, \qquad |\alpha| + |\beta| \geq k
\]
where $t \to (x^t,\xi^t)$ is the associated Hamilton flow. The case
$k=2$ has been considered in \cite{mmt0}; as proved there, the Hamilton
flow is bilipschitz and a wave packet parametrix can be constructed.
The case $k=0$ was considered in \cite{bt}; then the term $a^w(t,x,D)$ 
is  perturbative, and one may use the $a=0$ Hamilton flow in the above
condition. 

In the present article we need to deal with the case $k=1$. This
corresponds to a Hamilton flow which is close to the $a=0$ flow.
However, the term $a^w(t,x,D)$ is nonperturbative, and contributes a
time modulation factor along each packet.  Given these considerations,
we consider the following assumption on the symbol $a$:
\begin{equation}
\sup_{x,\xi} \int_{0}^1 |\partial_x^\alpha \partial_\xi^\beta a(t,x+2t\xi,\xi)| dt
\leq \epsilon c_{\alpha,\beta} \qquad |\alpha|+|\beta| \geq 1   
\label{l1a}\end{equation}

Let $(x_0,\xi_0)\in \R^{2n}$. To describe functions which are
localized in the phase space on the unit scale near $(x_0,\xi_0)$
we use the norm:
\[
H^{N,N}_{x^0,\xi^0}:=\{f: \langle D-\xi^0\rangle ^{N} f \in L^2,\ \langle x-x^0
\rangle^{N} f \in L^2 \}
\]

% If $f$ is localized in frequency at origin (neighborhood of size
% $\approx 1$) and highly localized in space at $x^0$ (neighborhood of
% size $\approx 1$), then $f \in H^{N,N}_{x^0}$ for any $N \in \N$.

We work with the lattice $\Z^n$ both in the physical and Fourier space.
 We consider a partition of unity in the physical space,
\[
\sum_{x_0 \in \Z^n} \phi_{x_0}=1  \qquad \phi_{x_0}(x) = \phi(x-x_0)
\]
where $\phi$ is a smooth bump function with compact support.
We  use a similar partition of unity on the Fourier side:
\[
\sum_{\xi_0 \in \Z^n} \varphi_{\xi_0} =1, \qquad \varphi_{\xi_0}(\xi) = \varphi(\xi-\xi_0)
\]
An arbitrary function $u$ admits an almost orthogonal 
decomposition
\[
u = \sum_{(x_0,\xi_0) \in \Z^{2n}} u_{x_0,\xi_0}, \qquad u_{x_0,\xi_0}
= \varphi_{\xi_0}(D) (\phi_{x_0} u)
\]
so that 
\begin{equation}
\sum_{(x_0,\xi_0) \in \Z^{2n}}
\|u_{x_0,\xi_0}\|_{H^{N,N}_{x^0,\xi^0}}^2 \lesssim \|u\|_{L^2}^2
\end{equation}
We remark that a continuous analog of the above discrete
decomposition can be obtained using the Bargman transform.

We first establish that the Hamilton flow is close to the Hamilton
flow with $a=0$:

\begin{l1}
Assume that \eqref{l1a} holds with a small
enough $\epsilon$. Then for each $(x^0,\xi^0) \in \R^{2n}$ and $t \in
[0,1]$ we have 
\[
|x^t - (x^0+2t\xi^0)| + | \xi^t -\xi^0| \lesssim \e
\]
\end{l1}

The proof is straightforward and is left for the reader; it is also
essentially contained in \cite{bt}.  This allows us to apply the main
result in \cite{mmt0}:

\begin{p1} \label{propsol} Assume that \eqref{l1a} holds with a small
  enough $\epsilon$.  Then for each $N \geq 0$ the solution of the
  homogeneous problem \eqref{generic} satisfies the following
  localization estimate:
\begin{equation} \label{locest} 
\| u(t) \|_{H^{N,N}_{x_0+2t\xi_0,\xi_0}} \les_N \| u_0 \|_{H^{N,N}_{x_0,\xi_0}}
\end{equation}
\end{p1}
We denote the evolution operator for \eqref{generic} by $S(t,s)$.
If the initial data is $u_0 = \delta_x$ then it has a decomposition of
the form
\begin{equation}
u_0 = \sum_{\xi_0 \in \Z^n} u_{\xi_0}(0), \qquad
\|u_{\xi_0}(0)\|_{H^{N,N}_{x,\xi^0}} \lesssim 1
\label{decd}\end{equation}
By \eqref{propsol}, at  time $1$ the corresponding solutions $u_{\xi_0}$ are
concentrated close to $x+2t\xi_0$, therefore they are spatially separated.
Hence we obtain the following pointwise decay:

\begin{c1} \label{cordecay}
The kernel $K(1,0)$ of $S(1,0)$ satisfies
\[
|K(1,x,0,y)| \lesssim 1
\]
The solution of the homogeneous equations \eqref{generic} satisfies
\[
\|S(1,0) u_0\|_{L^\infty} \lesssim \|u_0\|_{L^1}
\]
\end{c1}

If in addition the initial data is localized at some frequency
$\lambda$, say $u_0 = S_\lambda \delta_x$ then the decomposition in
\eqref{decd} is restricted to the range $|\xi_0| \approx \lambda$.
Then the corresponding solutions travel with speed $O(\lambda)$, and
we can obtain better pointwise decay away from the propagation region:

\begin{c1} \label{cordecay1}
The kernel $K_\lambda (1,0)$ of $S(1,0) S_\lambda$ satisfies
\begin{equation}
|K(1,x,0,y)| \lesssim (\lambda+|x-y|)^{-N}, \qquad |x-y| \not\approx \lambda 
\label{kernelbd1}\end{equation}
The kernel $K_\lambda (t,s)$ of $S(t,s) S_\lambda$ satisfies
\begin{equation} \label{kernelbd2}
|K(t,x,s,y)| \lesssim \lambda^{-N}, \qquad    |x-y| \approx \lambda, |t-s| \ll 1
\end{equation}

\end{c1}

The next result concerns localized energy estimates. 

\begin{c1} 
For each ball $B_r$ of radius $r \geq 1$ the solution $u$ to \eqref{generic} satisfies
\begin{equation} \label{lee}
\| S(t,0) S_\lambda u_0\|_{L^2(B_r)} \lesssim \lambda^{-\frac12}
r^\frac12 \|u_0\|_{L^2}
\end{equation}
\end{c1}

\begin{proof}
We consider the wave packet decomposition for $u = S(t,0) S_\lambda
u(0)$, 
\[
u = \sum_{(x_0,\xi_0) \in \Z^{2n}}^{|\xi_0| \approx \lambda} u_{x_0,\xi_0} 
\]
Let $\chi_r$ be a cutoff corresponding to $B_r$. Since $r \geq 1$ it
follows that the functions $\chi_r u_{x_0,\xi_0} $ are almost
orthogonal, therefore it suffices to prove the estimate for a single
packet. But a single packet is concentrated near a tube of spatial
size $1$ which travels with speed $O(\lambda)$. This tube intersects
the cylinder $[0,1] \times B_r$ over a time interval of length
$\lambda^{-1} r$. The conclusion easily follows.
\end{proof}

To obtain any results below the unit spatial scale we slightly strengthen 
the condition \eqref{l1a} by adding a weaker pointwise bound
\begin{equation}
|\partial_x^\alpha \partial_\xi^\beta a(t,x,\xi)|  \leq c_{\alpha \beta} \langle\xi\rangle^\frac12, \qquad \forall \ \alpha, \beta   
\label{l1b}\end{equation}
This will guarantee that on a unit spatial scale the flow in
\eqref{generic} is a small perturbation of the flat Schr\"odinger flow.
Then we have:

\begin{p1}
Assume that the conditions \eqref{l1a} and \eqref{l1b} hold. Then  

(i) For any $r > 0$ the solution $u$ to \eqref{generic} satisfies
the localized energy estimates
\begin{equation} \label{ler}
\| S(t,0) S_\lambda u(0)\|_{L^2(B_r)} \lesssim \lambda^{-\frac12}
r^\frac12 \|u_0\|_{L^2}
\end{equation}

(ii) For each $y,z$ with $|y-z| \approx \lambda$ we have
the square function bound
\begin{equation} \label{cordecay2}
\left \| \int_I S_\lambda S(t,s) S_\lambda (f(s) \delta_y) \right
\|_{L^2} \lesssim \lambda^{-1} \|f\|_{L^2}   
\end{equation}
\end{p1}

\begin{proof}
To prove this result it is convenient to replace the $L^2$ initial
data space by weighted $L^2$ spaces.

\begin{d1}
A weight $m: \R^{2n} \to \R^+$ is admissible if
\[
|m(x,\xi)/m(y,\eta)| \lesssim (1+|x-y|+|\xi-\eta|)^N
\]
for some real $N$.
\end{d1}

Correspondingly we define a weighted $L^2$ space
\[
\| u \|_{L^2(m)}^2 = \sum_{(x_0,\xi_0) \in \Z^{2n}} \| m(x_0,\xi_0)
u_{x_0,\xi_0}\|_{H^{N,N}_{x_0,\xi_0}}^2
\]
Given a weight $m_0$ at time $0$ we evolve it in time by
\[
m_t(x+2t\xi,\xi) = m_0(x,\xi)
\]
As a consequence of Proposition~\ref{propsol} we obtain

\begin{l1}
Assume that \eqref{l1a} holds with a small
enough $\epsilon$. Then 
\[
\| S(t,s)\|_{L^2(m_s) \to L^2(m_t)} \lesssim 1
\]
\end{l1}

Next we consider truncated solutions on a unit spatial scale. Given a
unit ball $B$ and an associated cutoff function $\chi$ we have the
following weighted local energy estimates:

\begin{l1} \label{llocal}
For any solution $u$ to \eqref{generic} we have
\[
\| \chi u\|_{L^2([0,1], L^2(\langle \xi \rangle^\frac12 m_t))}
+ \| (i \partial_t -\Delta)\chi u\|_{L^2([0,1], L^2(\langle \xi \rangle^{-\frac12} m_t))}
\lesssim \|u_0\|_{L^2(m)}
\]
\end{l1} 

\begin{proof}
We begin again with a wave packet decomposition of $u$,
\[
u = \sum_{(x_0,\xi_0) \in \Z^{2n}} u_{x_0,\xi_0} 
\]
The functions $ \chi u_{x_0,\xi_0}$ are almost orthogonal in $L^2$
therefore the bound for $\chi u$ follows.  On the other hand we have
\[
(i \partial_t -\Delta)\chi u = -2 \nabla \chi \nabla u - \Delta \chi u
+ \chi a^w(t,x,D)u
\]
The first two terms are estimated using the bound for $\chi u$.
For the last one we note that, by \eqref{l1b}, the operator
 $a^w$ preserves the $H^{N,N}$ spaces,
\[
\| a^w(t,x,D) u\|_{H^{N-1,N-1}_{x_0,\xi_0}} \lesssim \langle \xi_0
\rangle^\frac12 \| u\|_{  H^{N,N}_{x_0,\xi_0}} 
\]
Hence  we can use orthogonality again.
\end{proof}

Now we can conclude the proof of the Proposition.  For the local
energy estimate \eqref{lee} we first truncate $u$ to a unit scale. By
the above lemma with $m= (1+ \lambda^{-3} \langle \xi
\rangle^3)(1+\lambda^{3} \langle \xi \rangle^{-3})$ we obtain
\[
\lambda^{\frac12} \|\chi_1  u\|_{L^2(m')} + \lambda^{-\frac12}   \|(i
\partial_t -\Delta)\chi_1  u\|_{L^2(m')}  \lesssim \|u_0\|_{L^2(m)}
\]

where $m'=(1+ \lambda^{-2} \langle \xi \rangle^2)(1+\lambda^{2}
\langle \xi \rangle^{-2})$. It remains to show that
\[
\lambda^{\frac12} r^{-\frac12} \|\chi_r u\|_{L^2} \lesssim
\lambda^{\frac12} \|\chi_1 u\|_{L^2(m')} + \lambda^{-\frac12} \|(i
\partial_t -\Delta)\chi_1  u\|_{L^2(m')}
\]
Then we can localize the right hand side to the $\lambda^{-1}$ time
scale. On the $\lambda^{-1}$ time scale we can use the Duhamel formula
to further reduce the problem to a corresponding estimate for
solutions to the homogeneous constant
coefficient Schr\"odinger equation, namely:
\[
\lambda^{\frac12} r^{-\frac12} \|\chi_r e^{-it\D} u_0 \|_{L^2} \lesssim \| u_0 \|_{L^2(m')}
\]
After a dyadic frequency decomposition this becomes
\[
\lambda^{\frac12} r^{-\frac12} \|\chi_r e^{-it\D} S_\lambda u_0 \|_{L^2} \lesssim \| u_0 \|_{L^2}
\]
which is exactly the local energy estimate for the homogeneous constant
coefficient Schr\"odinger equation.

Consider now the square function bound. For $|t-s| \ll 1$ we can use
the kernel bound \eqref{kernelbd2}. Hence without any restriction in
generality we assume that $t \in I$, $s \in J$ where $I$, $J$ are
intervals of size $O(1)$ with $O(1)$ separation.
Choose  $t_0$ the center of the interval between $I$ and $J$.
We factor the estimate in two and prove the dual estimates
\[
\left \| \int_I S(t_0,s) S_\lambda (f(s) \delta_y) \right
\|_{L^2(m)} \lesssim \lambda^{-\frac12} \|f\|_{L^2}  
 \]
respectively
\[
\| (S_\lambda S(t,t_0) u)(z)\|_{L^2} \lesssim  \lambda^{-\frac12} \|u\|_{L^2(m)}
\]
where the flow invariant weight $m$ is given by 
\[
m(x,\xi) = (1+\lambda^{-1} |\xi\wedge(x-y)|)^{K} (1+\lambda^{-1} |\xi\wedge(x-z)|)^{-K}
\]
with $K$ large enough. These are dual bounds therefore it suffices 
to prove the second one. 

If $\chi$ is a smooth approximation of the characteristic function of
$B(z,1)$, then by (a slight modification of) Lemma~\ref{llocal} it remains to show that $v = \chi
S_\lambda u$ satisfies
\[
\lambda^\frac12 \| v(t,x)\|_{L^2(J)} \lesssim \lambda^{\frac12} \|
v\|_{L^2_t L^2(m \cdot m')} + \lambda^{-\frac12} \|(i
\partial_t -\Delta) v\|_{L^2_t L^2(m \cdot m')}
\]
where the additional weight 
$m'=(1+ \lambda^{-2} \langle \xi \rangle^2)(1+\lambda^{2} \langle \xi
\rangle^{-2})$ can be added due to the localization to frequency
$\lambda$.

This estimate can be localized to the $\lambda^{-1}$ timescale. In
addition, since $v$ has support in $B(z,1)$ we can freeze $x=z$ in $m$
and replace $m$ by
\[
\tilde{m}(\xi) = (1+\lambda^{-1} |\xi\wedge(y-z)|)^{K}
\]
Assuming $y-z = O(\lambda) e_1$ we get
\[
\tilde{m}(\xi) = (1+ |\xi'| + \lambda^{-1} |\xi_1|)^{K} 
\]
Then the $x'$ variable can be factored out and we are left with a
bound for the
one dimensional Schr\"odinger equation,
\[
 \| e^{-it\D}v_0 (\cdot,0)\|_{L^2(J)} \lesssim \lambda^{-\q} \|
v_0\|_{ L^2(m')}
\]
But this is exactly the one dimensional local energy estimate.

\end{proof} 
\section{The short time structure}
\label{scalel}

In this section we  consider a paradifferential approximation to
the magnetic Schr\"odinger equation \eqref{schra}. Precisely
given a dyadic frequency $\lambda$ we consider the evolution
\begin{equation}
i u_t - \Delta u + i (A_{< \sqrt \lambda} \nabla \tilde S_\lambda+ 
\tilde S_\lambda A_{< \sqrt \lambda} \nabla)
 u  = 0, \qquad u(0) = u_0
\label{schralambda}\end{equation}
where
\[
A_{< \sqrt \lambda} = S_{< \sqrt \lambda} A
\]
The multiplier $\tilde S_\lambda$ is added here for convenience.  It
guarantees that waves at frequencies away from $\lambda$ evolve
according to the constant coefficient Schr\"odinger flow, thereby
strictly confining the interesting part of the evolution to frequency
$\lambda$. In addition, the above expression is written in a
selfadjoint form, which guarantees that the corresponding evolution
operators $S(t,s)$ are $L^2$ isometries.

Later we will prove that on the time scale $\lambda^{-1}$ the
evolution of the $\lambda$ dyadic piece of a solution $u$ to
\eqref{schra} is well approximated by the evolution in
\eqref{schralambda}. Here we establish dispersive type estimates for
\eqref{schralambda}. Our main result concerning the flow in
\eqref{schralambda} is as follows:

\begin{p1} \label{pst} Let $u_\lambda$ be the solution to
  \eqref{schralambda} with initial data $u_{0,\lambda}$ localized at
  frequency $\lambda$.  Then for any interval $I$ of size less than
  $\lambda^{-1}$ the following estimates hold:

(i) the full Strichartz estimates
\begin{equation} \label{stri}
\|u_\lambda \|_{L^p(I,L^q)} \lessa \|u_{0,\lambda}\|_{L^2}
\end{equation}

(ii) the square function estimate
\begin{equation} \label{sf}
\|u_\lambda\|_{L^4_x(L^2_t(I))} \lessa \lambda^{-\frac14} \|u_{0,\lambda}\|_{L^2}
\end{equation}

(ii) the localized energy estimate: for any ball $B_r$ of radius $r >
0$ we have
\begin{equation} \label{lelambda}
\|u_\lambda\|_{L^2 (I \times B_r)} \lessa
r^\frac12 \lambda^{-\frac12}
 \|u_{0,\lambda}\|_{L^2}
\end{equation}
\end{p1}
The equation \eqref{schralambda} is $L^2$ well-posed, therefore we can
define the spaces $U^2_{A,\sqrt\lambda} L^2$, respectively
$V^2_{A,\sqrt\lambda} L^2$.

 As a consequence of the above 
proposition we obtain

\begin{c1}
  Assume that $A \in U^2_W H^1$. Then for any interval $I$ of length
  $\leq \l^{-1}$ and any function $u_\lambda$ localized at frequency
  $\lambda$ the following embeddings hold:
\begin{equation}
\| u_\lambda \|_{L^p(I,L^q)} \lessa \|u_{\lambda}\|_{U^2_{A,\sqrt{\lambda}} L^2}
\label{embarl}\end{equation}

\begin{equation}
\| u_\lambda \|_{L^4_x(L^2_t(I))}  \lessa \lambda^{-\frac14}
\|u_{\lambda}\|_{U^2_{A,\sqrt{\lambda}} L^2}
\label{embarlsf}\end{equation}

\begin{equation}
\| u_\lambda \|_{L^2 (I \times B_r)} \lessa r^\frac12 \lambda^{-\frac12}
\|u_{\lambda}\|_{U^2_{A,\sqrt{\lambda}} L^2}
\label{embarlle}\end{equation}

\end{c1}
The first step in the proof of Proposition~\ref{pst} is to establish a
wave packet parametrix and a wave packet decomposition for solutions
to \eqref{schralambda} on the $\lambda^{-1}$ time scale. This is done
by rescaling starting from the results in the previous section.

We begin by writing in the Weyl  calculus
\[
i (A_{< \sqrt \lambda} \nabla \tilde S_\lambda+ 
\tilde S_\lambda A_{< \sqrt \lambda} \nabla) = a^w(t,x,D)
\]
Then the symbol $a(t,x,\xi)$ can be expressed as a principal term 
plus an error,
\[
a(t,x,\xi) = a_0(t,x,\xi) + a_r(t,x,\xi) 
\]
where the principal part $a_0$ is given by
\[
a_0(t,x,\xi) = -2 i A_{< \sqrt \lambda}(t,x) \cdot \xi \tilde s_\lambda(\xi)
\]

By Sobolev embeddings we have the following pointwise bound 
for the truncated magnetic potential:
\[
| \partial_x^\alpha A_{< \sqrt \lambda}(t,x)| \leq c_\alpha
\lambda^\frac{1+|\alpha|}2 \|A(t)\|_{H^1}
\]
This yields 
\begin{equation}
| \partial_x^\alpha \partial_\xi^\beta a_0(t,x,\xi)| \lessa c_{\alpha \beta}
\lambda^{\frac{3+|\alpha|}2-|\beta|}
\label{a0bd}\end{equation}
In addition, by the Weyl calculus it follows that
$a_r$ is also localized at frequency $\lambda$ and satisfies
\begin{equation}
| \partial_x^\alpha \partial_\xi^\beta a_r(t,x,\xi)| \lessa c_{\alpha \beta}
\lambda^{\frac{1+|\alpha|}2-|\beta|}
\label{arbd}\end{equation}
This brings us to our main integral bound for the symbol $a$, namely

\begin{l1}
Assume that $A \in U^2_W H^1$ with $\text{div } A = 0$. Then the above
symbol $a$ satisfies
 \begin{equation}
\sup_{x,\xi} \int_{0}^{T}  |\partial_x^\alpha \partial_\xi^\beta a(t,x+2t\xi,\xi)| dt
\lessa
\left\{ \begin{array}{cc} c_{\beta} (T\lambda)^\frac12
    \lambda^{-|\beta|} \log \lambda & \alpha = 0
 \cr c_{\alpha\beta} (T\lambda)^\frac12 \lambda^{\frac{|\alpha|}2-|\beta|},
& |\alpha| \geq 1.
\end{array}\right.
\label{l1aa}\end{equation}
\label{lsyma}\end{l1}

\begin{proof}
The bound for $a_r$ follows directly from \eqref{arbd}, therefore it
remains to consider $a_0$. Furthermore, it suffices to consider the
case $|\alpha|=0,1$, $\beta = 0$. Then we need to prove the bounds
\[
 \int_{0}^{T}  |A_{< \sqrt \lambda}(t,x+2t\xi)|  dt
\lesssim T^\frac12 \lambda^{-\frac12} \ln \lambda  
\|\nabla A\|_{U^2_W L^2}, \qquad |\xi| \approx \lambda
\]
respectively
\[
 \int_{0}^{T}  |\nabla A_{< \sqrt \lambda}(t,x+2t\xi)|  dt
\lesssim T^\frac12 \|\nabla A\|_{U^2_W L^2}, \qquad |\xi| \approx \lambda
\]
These in turn follow by dyadic summation from
\begin{equation}
\sup_{x,\xi} \int_{0}^{T}  |(S_\mu B)(t,x+2t\xi)|  dt
\lesssim T^\frac12 \mu \lambda^{-\frac12} \|B\|_{U^2_W L^2}, \qquad |\xi| \approx \lambda
\label{kiu}\end{equation}

The line $y = x+2t\xi$ moves through a unit spatial cube in a time
$\lambda^{-1}$. But, due to the finite speed of propagation for the wave
equation, see \eqref{wcloc}, the contributions from different spatial unit cubes are
square summable. Hence by Cauchy-Schwartz it suffices to
prove the above bound for $T \leq \lambda^{-1}$.
By \eqref{pkjh}, for $ T \leq \lambda^{-1}$ we have 
\[
\|B\|_{U^2_W (0,T;L^2)} \approx \|B\|_{U^2 (0,T;L^2)},  
\]
therefore it is enough to prove: 
\[
 \int_{0}^{T}  |(S_\mu B)(t,x+2t\xi)|  dt
\lesssim (\l T)^\frac12 \mu \l^{-1} \|B\|_{U^2 L^2}
\]
It suffices to prove the bound when $B$ is an $U^2 L^2$ atom. 
By Cauchy-Schwartz it suffices to consider a single step,
which corresponds to a time independent $B$. Then the last bound 
can be rewritten in the form
\[
\int_L  |(S_\mu B)(x)| ds \lesssim \mu |L|^\frac12 \|B\|_{L^2}
\]
where $L$ is an arbitrary line segment in $\R^3$. We can set $\mu
= 1$ by rescaling. In coordinates $x = (x_1,x')$ suppose $L$ is
contained in $\{x'=0\}$. Then we use Sobolev embeddings in $x'$ and 
Cauchy-Schwartz with respect to $x_1$.
\end{proof}

Next we consider the rescaling that preserves the flat Schr\"odinger
flow and takes the time scale $\lambda^{-1}$ to $1$, namely
\[
 v_\lambda(x,t)=u\left(\frac{x}{\sqrt{\l}},\frac{t}{\l}\right)
\]
If $u$ solves \eqref{schralambda} then for $v_\lambda$ we obtain the following
equation: 
\begin{equation} \label{rescaled} i v_t - \Delta v +
  \lambda^{-1}a^w\left(\frac{t}{\l},\frac{x}{\sqrt{\l}},D \sqrt \lambda\right) v =    0
\end{equation} 
However, this is not sufficient, we need to repeat the same procedure
for shorter time scales. Precisely, for each $ \lambda < \mu <
\lambda^2$ we can rescale the $\mu^{-1}$ time scale to the unit scale
by setting
\[
 v(x,t)=u\left(\frac{x}{\sqrt{\mu}},\frac{t}{\mu}\right)
\]
Then for $v$ we obtain the equation 
\begin{equation} \label{rescaled2} i v_t - \Delta v + a_\mu^w\left(
    t,x,D) \right) v = 0, \quad a_\mu(t,x,\xi) =
  \mu^{-1}a^w\left(\frac{t}{\mu},\frac{x}{\sqrt{\mu}},\xi \sqrt \mu
  \right)
\end{equation}

Rescaling the bounds \eqref{a0bd}, \eqref{arbd} and \eqref{l1aa} it
follows that this rescaled equation belongs to the class studied in
the previous section:

\begin{l1} \label{symbolcheck} For $ \e^{-1} \lambda \leq \mu \leq \lambda^2$
  and $\epsilon$ small enough the symbol $a_\mu$ satisfies \eqref{l1a}
  and \eqref{l1b} on the time interval $[0,1]$.
\end{l1}

This allows us to apply the results in the previous section to the
evolution~\eqref{schralambda}.  Rescaling the result in Corollary
\ref{cordecay} we obtain short time pointwise bounds for the solution
to \eqref{schralambda}:

\begin{l1}
The solution of \eqref{schralambda} has the  pointwise decay
\begin{equation} \label{pointws}
\| u(t) \|_{L^{\infty}} \les |t-s|^{-\frac{n}{2}} \| u(s) \|_{L^1},
\qquad |t-s| \les \epsilon \l^{-1}
\end{equation}
\end{l1}

\begin{proof}
W.a.r.g we can take $s = 0$.
  If $  \lambda^{-2} \leq t \leq \e \lambda^{-1}$ then this follows
  directly from Corollary~\ref{cordecay} applied to the equation
  \eqref{rescaled2} with $ \mu^{-1}=t$.  

  The case $t <  \lambda^{-2}$ needs to be considered separately.
For such $t$ we split the evolution in two parts,
\[
S(t,0) = S(t,0) \tilde{\tilde S}_{\l} + S(t,0) (1-\tilde{\tilde S}_{\l})
\]
The second part evolves according to the constant coefficient
Schr\"odinger flow, hence it is easy to estimate. For the first part
we use the rescaled parametrix in the previous section corresponding
to $\mu = \lambda^{-2}$. The solution $S(t,0) \tilde{\tilde S}_{\l} \delta_x$ 
consists of a single packet on the $\lambda^{-1}$ spatial scale
which does not move up to time $\lambda^{-2}$. 
Hence we obtain
\begin{equation}
|S(t,0) \tilde{\tilde S}_{\l} \delta_x| \lesssim \lambda^{3}
(1+\lambda |x-y|)^{-N}, \qquad |t| \leq \lambda^{-2}
\label{lowtbd}\end{equation}
which concludes the proof.
\end{proof}

By \cite{MR1646048}, the Strichartz estimates in \eqref{stri} are a
direct consequence of \eqref{pointws}. We continue with a 
decay bound away from the propagation region:

\begin{l1} 
If $|t-s| \leq \e \lambda^{-1}$
then the kernel of $S(t,s) S_\lambda$ satisfies
\begin{equation}
|K(t,x,s,y)| \lesssim \lambda^3 (1+  \lambda|x-y|+ \lambda^2|t-s| )^{-N}
\label{fardecay}\end{equation}
whenever
\[
|t-s|+\lambda^{-2}
\gg \lambda^{-1}|x-y|\ \ \  \ \text{or} \ \ \ \ \lambda^{-1}|x-y| +
\lambda^2 \gg |t-s| 
\]
\end{l1}

\begin{proof} W.a.r.g we can take $s = 0$.  If $|t-s| \leq
  \lambda^{-2}$ then we use \eqref{lowtbd}.  On the other hand if $
  \lambda^{-2} \leq |t-s| \leq \e \lambda^{-1}$ and then we rescale
  \eqref{kernelbd1} applied to \eqref{rescaled2} with $\mu = t^{-1}$.

Since the input is localized at frequency $\lambda$, it
follows that waves need exactly a time $\approx \lambda^{-1}|x-y|
+\lambda^{-2}$ to travel from $x$ to $y$.

\end{proof}

Next we consider pointwise 
square function bounds:

\begin{l1}
  The evolution $S(t,s)$ associated to
  \eqref{schralambda} has the pointwise square function decay 
\begin{equation} \label{sqfun}
  \left\| \int_I S(t,s) S_\lambda (f(s) \delta_y)ds
    (x) \right\|_{L^2_t(I)} \les |x-y|^{-1} \| f \|_{L^2_t(I)}, \quad
  |I| \les \epsilon \l^{-1} 
\end{equation}
\end{l1}
\begin{proof}
If $|x-y| \gtrsim 1$ then we can use directly \eqref{fardecay}. If
$|x-y| \ll 1$ then we split the integral in two parts. If $|t-s| \gg
\lambda^{-1} |x-y|$ then we can still use  \eqref{fardecay}. On the
other hand if $|t-s| \lesssim \lambda^{-1} |x-y|$ then we rescale
 \eqref{cordecay2} applied to \eqref{rescaled2} with $\mu = \lambda |x-y|^{-1}$.
\end{proof}

We continue with the proof of \eqref{sf}. By the $TT^*$ argument we
need to prove the bound
\begin{equation}
\left \|  \int_I S(t,s) S_\lambda f(s) ds
\right\|_{L^4_x L^2_t} \lessa \lambda^{-\frac12} \| f\|_{L^\frac43_x L^2_t} 
\label{sf1}\end{equation}
For this we use Stein's complex interpolation theorem. Define the
holomorphic family of operators
\[
T_z f (t) =  \int_I  z (t-s)_{+}^{z-1} S(t,s) S_\lambda f(s) ds
\]
Then we need to show that 
\[
\| T_1 f\|_{L^4_x L^2_t} \lessa  \lambda^{-\frac12} \| f\|_{L^\frac43_x L^2_t} 
\]
This follows by interpolation from
\begin{equation}
\| T_z f\|_{L^2} \lesssim \|f\|_{L^2}, \qquad \Re z = 0 
\label{sf2}\end{equation}
and
\begin{equation}
\| T_z f\|_{L^\infty_x L^2_t} \lessa \lambda^{-1} \|f\|_{L^1_x L^2_t}, \qquad \Re z = 2 
\label{sf3}\end{equation}

For \eqref{sf2} we write
\[
S(0,t) T_z f(t) =  \int_I  z (t-s)_{+}^{z-1}
S_\lambda S(0,s) f(s) ds
\]
Since $S(t,s)$ are $L^2$ isometries it suffices to
prove that 
\[
\|  \int_I  z (t-s)_{+}^{z-1} f(s) ds\|_{L^2} \lesssim \|f\|_{L^2}
\]
which is straightforward by Plancherel's theorem since the Fourier transform of $z
t_{+}^{z-1}$ is $\Gamma(z+1) (\tau+i0)^{-z}$ which is bounded.

On the other hand the bound \eqref{sf3} is equivalent to
\[
\| T_z (f \delta_y) (x)\|_{L^2_t} \lessa  \lambda^{-1} \|f\|_{L^2_t} 
\]
which we can rewrite in the form
\[
\left \|  \int_I  (t-s)^{1+i\sigma} S_{\sqrt{\lambda}}(t,s) S_\lambda (f(s) \delta_y)ds
    (x) \right\|_{L^2_t(I)} \lessa \lambda^{-1} \| f \|_{L^2_t(I)}
\]
Restricting $t-s$ to the range $|t-s| \lesssim \lambda^{-1}|x-y| +\lambda^{-2}$ this
is a consequence of \eqref{sqfun}. On the other hand for larger $t-s$ we
can use directly the pointwise bound \eqref{fardecay}.

The last step of the proof of Proposition~\ref{pst} is the localized
energy estimate \eqref{lelambda}. This follows directly by rescaling
from \eqref{ler}  applied to the equation \eqref{rescaled}.

\section{ Short range bilinear and trilinear estimates}

We first consider $L^2$ bilinear product estimates where one 
factor solves the wave equation and the other solves the 
Schr\"odinger equation.

\begin{p1}
Assume that $A \in U^2_W H^1$, $1 \leq \mu \lesssim \lambda$ and
$|I| \leq \lambda^{-1}$. Then the following bilinear 
$L^2$ estimates hold:
\begin{equation} \label{tri3}
\| S_\mu (B_\l u_\l) \|_{L^2(I \times \R^3)} \lessa \mu^{\frac12}  
\| B_\l \|_{U^2_W L^2} \| u_\l \|_{U^2_{A,\sqrt\lambda} L^2},
\end{equation} 
\begin{equation} \label{tri4}
\| B_\l u_\mu \|_{L^2(I \times \R^3)} \lessa \mu^{\frac12}  
\| B_\l \|_{U^2_W L^2} \| u_\mu \|_{U^2_{A,\sqrt\mu} L^2} 
\end{equation} 
\begin{equation} \label{tri2}
\| B_\mu  u_\lambda \|_{L^2(I \times \R^3)} 
\lessa \mu \l^{-\frac12}  \| B_\mu \|_{U^2_W L^2} 
\| v_\lambda \|_{U^2_{A,\sqrt\lambda}
  L^2},
\end{equation}
\end{p1}

We remark that the constants in \eqref{tri2} are optimal,
and in effect  as a consequence of the results in the
last section \eqref{tri2} can be extended almost up to time $1$.
On the other hand the constants in \eqref{tri3}, \eqref{tri4} 
are not optimal, but this is not so important because this corresponds
to the non-resonant case in the trilinear estimates.

\begin{proof}
  For \eqref{tri3} it suffices to use Bernstein's inequality and the
  Strichartz estimates,
\[
\begin{split}
\| S_\mu (B_\l u_\l) \|_{L^2}
\lesssim &\ \mu^\frac12 \|B_\l u_\l \|_{L^2_t L^\frac32_x}
\\ \lesssim &\ \mu^\frac12 \|B_\l\|_{L^\infty L^2} \| u_\l \|_{L^2_t L^6_x}
\\ \lessa &\ \mu^{\frac12}  
\| B_\l \|_{U^2_W L^2} \| u_\l \|_{U^2_{A,\sqrt\lambda} L^2}
\end{split}
\]
A similar argument applies for \eqref{tri4}.

It remains to prove \eqref{tri2}.  By \eqref{pkjh} we 
replace the $U^2_W L^2$
space by the $U^2 L^2$ space on a short time scale.
Hence we can rewrite \eqref{tri2} in the form
\begin{equation}
\| B_\mu  u_\lambda \|_{L^2(I \times \R^3)} 
\lessa \mu \l^{-\frac12}  \| B_\mu \|_{U^2 L^2} 
\| v_\lambda \|_{U^2_{A,\sqrt\lambda}
  L^2}
\label{tri2a}\end{equation}
Due to the atomic structure of the $U^2$ spaces 
it suffices to prove the above bound in the special case when
both $ B_\mu$  and  $u_\lambda$ solve the corresponding 
homogeneous equations $\partial_t B_\mu = 0$, respectively 
\eqref{schralambda}.

We consider a partition of unit on the $\mu^{-1}$ scale
\[
1 =  \sum_{x_0 \in \mu^{-1} \Z^3} \phi_{x_0}^2(x)
\]
and use the localized energy estimates \eqref{embarlle} for
$u_\lambda$  with $r = \mu^{-1}$:
\[
\begin{split}
\|  B_\mu  u_\lambda \|_{L^2(I \times \R^3)}^2&\  \approx 
 \sum_{x_0 \in \mu^{-1} \Z^3} \| \phi_{x_0}^2 B_\mu  u_\lambda\|^2_{L^2}
\\ &\ \lesssim  \sum_{x_0 \in \mu^{-1} \Z^3} \|  \phi_{x_0}
B_\mu\|_{L^\infty}^2  \|  \phi_{x_0} u_\lambda\|^2_{L^2}
\\ &\ \lesssim_A  \lambda^{-1} \mu^{-1} 
\|  u_\lambda\|^2_{U^2_{A,\sqrt\lambda} L^2}
 \sum_{x_0 \in \mu^{-1} \Z^3} \|  \phi_{x_0}
B_\mu\|_{L^\infty}^2
\\ &\ \lesssim  \lambda^{-1} \mu^2  \| B_\mu \|_{U^2 L^2} 
\| v_\lambda \|_{U^2_{A,\sqrt\lambda}
  L^2}
\end{split}
\]
\end{proof}

Next we turn our attention to trilinear estimates. We begin with the
easier case of three $U^2$ type spaces

\begin{p1} \label{u2tri}
a) If $|I| \leq \l^{-1}$, $\mu \les \l$ and $B_\mu$,
  $u_\lambda$, $v_\lambda$ are
  localized at frequency $\mu$, $\lambda$, respectively $\lambda$ then
  \begin{equation} \label{llmu}
\!\!\!\!  \left| \int_{I} \!\! \int_{\R^3}\! B_\mu u_\lambda \bar{v}_\lambda dx dt
  \right| \lessa  \! \frac{\min{(\mu,\l^\q)}}{\l} \| B_\mu\|_{U^2_W L^2}
  \|u_\lambda \|_{U^2_{A,\sqrt\lambda} L^2} \|
  v_\lambda\|_{U^2_{A,\sqrt\lambda} L^2}
  \end{equation}

  b) If $|I|=\l^{-1}$, $\mu \ll \l$ and $B_\lambda$, $u_\mu$,
  $v_\lambda$ are localized at frequency $\lambda$, $\mu$,
  respectively $\lambda$ then
  \begin{equation} \label{lmlu}
  \left| \int_{I} \int_{\R^3} B_\l u_\mu \bar{v}_\l dx dt \right| \les_A
  \mu^\q \l^{-1} \| B_\l \|_{U^2_W L^2} \|v_\mu
  \|_{U^2_{A,\sqrt\mu} L^2} \| w_\l \|_{U^2_{A,\sqrt\lambda} L^2}
  \end{equation}
\end{p1}

\begin{proof}
a) If $\mu < \lambda^\frac12$ then the conclusion follows 
directly from \eqref{tri2} since 
\[
 \left| \int_{I} \int_{\R^3} B_\mu u_\lambda \bar{v}_\lambda dx dt
  \right| \lesssim |I|^{\frac12} \| B_\mu u_\lambda\|_{L^2} \|
    v_\lambda\|_{L^\infty L^2}
\]

If $\mu > \lambda^\frac12$ then we use \eqref{pkjh} to replace 
$U^2_W L^2$ by $U^2 L^2 \subset L^2_x L^\infty_t$.
Then we  estimate
\[
 \left| \int_{I} \int_{\R^3} B_\mu u_\lambda \bar{v}_\lambda dx dt
  \right| \lesssim_A \|B_\mu\|_{ L^2_x L^\infty_t} \| u_\lambda\|_{L^4_x
    L^2_t}  \| v_\lambda\|_{L^4_x
    L^2_t} 
\]
and use the square function bounds \eqref{sf} for the last two factors.

b)  In the Fourier space we obtain nontrivial contributions
when either all three time frequencies are $\ll \lambda^2$ or when at
least two of them are $\gtrsim \lambda^2$. More precisely, 
using smooth time multiplier cutoffs we can
write
\[
\begin{split}
\int_{I}\int_{\R^3}\!\!  B_\l u_\mu \bar{v}_\l & dx dt = 
\int_{I} \int_{\R^3} \chi_{\{|D_t| > \lambda^2/32\}} B_\l u_\mu
\bar{v}_\l dx dt
\\ &\ + \int_{I} \int_{\R^3} \chi_{\{|D_t| < \lambda^2/32\}}  B_\l \chi_{\{|D_t| > \lambda^2/32\}} u_\mu \bar{v}_\l dx dt
 \\ &\ + \int_{I} \int_{\R^3} \chi_{\{|D_t| < \lambda^2/32\}}B_\l \chi_{\{|D_t| <
  \lambda^2/32\}} u_\mu \chi_{\{|D_t| < \lambda^2/16\}} \bar{v}_\l dx dt
\end{split}
\]

Since the wave equation has constant coefficients,
for the first term we can bound the first factor in $L^2$,
\[
\|  \chi_{\{|D_t| > \lambda^2/32\}} B_\l\|_{L^2} \lesssim \lambda^{-1}
\|B_\lambda\|_{U^2_W L^2}
\]
On the other hand for the remaining product we  use the energy estimate for $v_\lambda$
and the $L^2 L^\infty$ bound for $u_\mu$.

We argue in a similar manner for the other two terms. The bilinear expressions
\[
\chi_{\{|D_t| < \lambda^2/32\}}B_\l \chi_{\{|D_t| <
  \lambda^2/32\}} u_\mu, \qquad S_\mu(\chi_{\{|D_t| < \lambda^2/32\}}  B_\l   \bar{v}_\l)
\]
can be estimated in $L^2$ using \eqref{tri4}, respectively
\eqref{tri3}. Hence it remains to bound in $L^2$ the high modulation
factors:

\begin{l1}
We have
\[
\| \chi_{\{|D_t| < \lambda^2/16\}} {v}_\l\|_{L^2} \lesssim_A \lambda^{-1}
\| v_\lambda \|_{U^2_{A,\sqrt{\lambda}} L^2}
\]
respectively
\[
\| \chi_{\{|D_t| > \lambda^2/32\}} u_\mu\|_{L^2} \lesssim_A \lambda^{-1}
\| u_\mu \|_{U^2_{A,\sqrt{\mu}} L^2}
\]
\end{l1}
\begin{proof}
In the case $A=0$ both bounds are trivial, the difficulty is to
accommodate the unbounded term involving $A$. We consider
the first bound only, as the argument for the second is similar.

Without any restriction in generality we can take $v_\lambda$ to be an
$U^2$ atom.  The kernels of the operators $\chi_{\{|D_t| <
  \lambda^2/16\}}$, respectively $\chi_{\{|D_t| > \lambda^2/32\}}$
decay rapidly on the $\lambda^{-2}$ time scale. 
Then it suffices to prove the estimate in two cases:

(i) $v_\lambda$ is supported in a $\l^{-2}$ time interval 
(this corresponds to steps of length $\l^{-2}$ and shorter).
Then the bound follows directly from the energy estimates 
and Holder's inequality.

(ii) $v_\lambda$ solves the homogeneous
equation \eqref{schralambda} on the time interval $I$ with  $
\lambda^{-2} \leq |I| \leq \l^{-1}$ (this
corresponds to steps of length $\l^{-2}$ and longer).  Then we can use
the bound \eqref{tri2} to estimate
  \[
  \| (i\partial_t -\Delta) v_\lambda\|_{L^2(I)} = \| A_{< \sqrt \lambda} \nabla \tilde
  S_\lambda v_\lambda\|_{L^2(I)} \lesssim_A \ln \lambda \ \lambda^\frac12
  \|v_\lambda\|_{L^\infty L^2}
  \]
Hence with $I = [t_0,t_1]$ we can write
\[
 (i\partial_t -\Delta) (\chi_I v_\lambda) = i v_\lambda(t_0) \delta_{t=t_0} -
 i v_\lambda(t_1) \delta_{t=t_1} + f_\lambda
\]

where $\chi_I$ is the characteristic function of $I$ and

\[
\| f_\lambda \|_{L^2} \lesssim_A \ln \lambda \ \lambda^\frac12
  \|v_\lambda\|_{L^\infty L^2}
\]
Hence working with the constant coefficient Schr\"odinger equation we 
obtain 
\[
\begin{split}
\| \chi_{\{|D_t| < \lambda^2/2\}} \chi_I {v}_\l\|_{L^2} \lesssim_A &
\lambda^{-1} (\|v_\lambda(t_0)\|_{L^2} + \|v_\lambda(t_1)\|_{L^2})
+\lambda^{-2} \|f\|_{L^2} 
\\ \lessa & \lambda^{-1}
  \|v_\lambda\|_{L^\infty L^2}
\end{split}
\]
which is exactly what we need.

\end{proof}
\end{proof}

Finally we turn our attention to the case when one of the three
$U^2$ spaces is replaced by a $V^2$ space:

\begin{p1} a) Let  $|I| \leq \l^{-1}$, $\mu \les \l$ and $B_\mu$,
  $u_\lambda$, $v_\lambda$ 
  localized at frequency $\mu$, $\lambda$, respectively $\lambda$.
Then
  \begin{equation} \label{vllma}
 \!\!\!\!  \left| \int_{I}  \int_{\R^3}  B_\mu  u_\lambda  \bar{v}_\lambda dx
    dt \right|  \lessa 
 \frac{(\l |I|)^\q \mu}{\lambda}   
\| B_\mu\|_{U^2_W    L^2} 
\|u_\lambda \|_{U^2_{A,\sqrt\lambda} L^2} 
\| v_\lambda\|_{V^2_{A,\sqrt\lambda} L^2}
  \end{equation}
If in addition $\lambda^\q \ll \mu \lesssim \lambda$ then
  \begin{equation} \label{llma}
  \left| \int_{I} \int_{\R^3}  B_\mu  u_\lambda  \bar{v}_\lambda dx
    dt \right| \lessa 
\frac{ \ln\left(\frac{\mu}{\sqrt\lambda}\right)}{\sqrt \lambda} 
\| B_\mu\|_{U^2_W    L^2} 
\|u_\lambda \|_{U^2_{A,\sqrt\lambda} L^2} 
\| v_\lambda\|_{V^2_{A,\sqrt\lambda} L^2}
  \end{equation}

  b) If $|I|=\l^{-1}$, $\mu \ll \l$  and $B_\lambda$, $u_\mu$, $v_\lambda$ are
  localized at frequency $\lambda$, $\mu$, respectively $\lambda$ then
  \begin{equation} \label{lmla}
  \left| \int_{I} \int_{\R^3} B_\l u_\mu \bar{v}_\l dx dt \right| \les_A
  \frac{\mu^\q  \ln{\lambda}}{\l} \| B_\l \|_{U^2_W L^2} \|u_\mu
  \|_{U^2_{A,\sqrt\mu} L^2} \| v_\l \|_{V^2_{A,\sqrt\lambda} L^2}
  \end{equation}
\end{p1}

\begin{proof}
a) Using the bilinear $L^2$ bound \eqref{tri2} for the
product of the first two factors we obtain
\[
 \left| \int_{I}  \int_{\R^3}  B_\mu  u_\lambda  \bar{v}_\lambda dx
    dt \right|  \lessa 
 \frac{(\l |I|)^\q \mu}{\lambda}   
\| B_\mu\|_{U^2_W    L^2} 
\|u_\lambda \|_{U^2_{A,\sqrt\lambda} L^2} 
\| v_\lambda\|_{L^\infty L^2}
\]
Then \eqref{vllma} follows due to the trivial embedding
$V^2_{A,\sqrt\lambda} L^2 \subset L^\infty_t L^2_x$.

On the other hand the LHS of \eqref{llma} can be estimated either as
above or as in \eqref{llmu}.  Then \eqref{llma} follows from the
decomposition
\[
V^2_{A,\sqrt\lambda} L^2 \subset
\ln\left(\frac{\mu}{\sqrt\lambda}\right)   U^2_{A,\sqrt\lambda} L^2 +
\left(\frac{\mu}{\sqrt\lambda}\right)^{-1} L^2
\]
We can factor out the \eqref{schralambda} flow by pulling functions
back to time $0$ along the flow. Then the above relation becomes
\[
V^2 L^2 \subset \ln \sigma \ U^2 L^2 + \sigma^{-1} \ L^\infty L^2, \qquad \sigma
\gg 1
\]
This in turn is true due to Lemma~\ref{principle}.

b) This follows from a similar argument to the one above and by using
\eqref{tri4} and \eqref{lmlu}.

\end{proof}

\section{The short time paradifferential calculus}

In this section we prove that, given a dyadic frequency $\lambda$, the
evolution of the $\lambda$ dyadic piece of a solution $u$ to
\eqref{schra} is well approximated by the evolution of the
paradifferential equation \eqref{schralambda} on time intervals of
size $\lambda^{-1}$. We also introduce different paradifferential
truncations
\begin{equation}
i u_t - \Delta u + i (A_{<\nu} \nabla \tilde S_\lambda+ 
\tilde S_\lambda A_{< \nu} \nabla)
 u  = 0, \qquad u(0) = u_0
\label{schranu}\end{equation}
and show that they all generate equivalent spaces.  The spaces
associated to \eqref{schranu} are denoted by $U^2_{A,\nu,\lambda}
L^2$, $U^2_{A,\nu,\lambda} L^2$, respectively $U^2_{A,\nu,\lambda}
L^2$. We refer to the above evolution as the $(A_{<\nu},\lambda)$
flow. 

A special case of the above equation is when $A_{<\nu}$ is replaced 
by $A_{\ll \lambda}$. We refer to that as the $(A_{\ll \l},\l)$ flow.
By a slight abuse of notation we denote the corresponding spaces 
by $U^2_{A,\lambda,\lambda} L^2$, etc.

\begin{p1} \label{paal}
a) For any interval $I$ with $|I| \leq \lambda^{-1}$  the solution $u$
to \eqref{schra} satisfies
\begin{equation}
\| S_\lambda u\|_{U^2_{A,\sqrt\lambda}(I;L^2)} 
 \lesssim_A \| u_0\|_{L^2}
\label{sla} \end{equation}
b) In addition, for any $\sqrt\lambda < \nu \ll \lambda$ we have
\begin{equation}
\|  u\|_{U^2_{A,\sqrt\lambda}(I;L^2)} \approx_A    \| u\|_{U^2_{A,\nu,\l}(I;L^2)}
\label{slec}\end{equation}
\end{p1}

The bound \eqref{sla} transfers easily to $U^2$ spaces:
\begin{c1}
 For any interval $I$ with $|I| \leq \lambda^{-1}$ we have
\begin{equation}
\| S_\lambda  u\|_{U^2_{A,\sqrt\lambda}(I;L^2)} \lessa    \| u\|_{U^2_{A}(I;L^2)}
\label{sluc}\end{equation}
\end{c1}

Combining this with \eqref{embarl} we immediately obtain the bounds
\eqref{seu2} in part (i) of Theorem~\ref{ta}:

\begin{c1} \label{plla}
The solution $u$ to the homogeneous magnetic Schr\"odinger equation 
\eqref{schra}  satisfies the Strichartz estimates \eqref{seu2}.
\end{c1}

The rest of the section is dedicated to the proof of Proposition \ref{paal}.

\begin{proof}
a) From the equation \eqref{schra} we obtain the following equation for 
$S_\lambda u$,
\[
\left(i \partial_ t - \Delta + i  A_{< \sqrt \lambda} \nabla \tilde
  S_\lambda   +  i\tilde S_\lambda   A_{< \sqrt \lambda} \nabla  \right) S_\lambda u =
f_\lambda
\]
where 
\[
f_\lambda = S_\lambda ( 2i A_{> \sqrt \lambda} \nabla u + A^2 u) + i [
S_\lambda , A_{< \sqrt \lambda} ] \tilde
  S_\lambda \nabla u 
\]

Then we have 
\[
\| S_\lambda u\|_{U^2_{A,\sqrt{\lambda}}(I;L^2)} \lesssim 
\|u_0\|_{L^2} + \| f_\l \|_{DU^2_{A,\sqrt{\lambda}}(I;L^2)} 
\]
The estimate \eqref{sla} follows if we establish that the
inhomogeneous terms $f_\lambda$ are uniformly small,
\[
\| f_\l \|_{DU^2_{A,\sqrt{\lambda}}(I;L^2)} \lessa (\lambda |I|)^\delta 
\sup_{\lambda'}  \| S_{\lambda'} u\|_{U^2_{A,\sqrt{\lambda'}}(I;L^2)}
\]
We consider the terms in $f_\lambda$. For the first term by duality we
need to prove that
\[
\left| \int_I \int_{\R^3}  A_{> \sqrt \lambda} \nabla u   S_\lambda \bar v dx
dt \right| \lessa   (\lambda |I|)^\delta  \|v\|_{V^2_{A,\sqrt{\lambda}} L^2}
\sup_{\lambda'}  \| S_{\lambda'} u\|_{U^2_{A,\sqrt{\lambda'}}(I;L^2)}
\]
We take a Littlewood-Paley decomposition of the first two factors and
estimate each dyadic piece. There are several cases to consider for
the integrand:

a) $S_\mu A \nabla S_\lambda u S_\lambda \bar v$ with $\sqrt\l \leq
\mu \leq \lambda$. Using \eqref{vllma} yields a constant
\[
|I|^\frac12 \sqrt \lambda = (\lambda |I|)^\frac14 ( \mu^2 |I|)^\frac14 
\left(\frac{\l}{\mu^2}\right)^\frac14
\]
which is favorable if $|I| \leq \mu^{-2}$. On the other hand using
\eqref{llma} yields a constant
\[
 \ln\left(\frac{\mu}{\sqrt\lambda}\right) \mu^{-1} {\sqrt \lambda} 
=  \ln\left(\frac{\mu}{\sqrt\lambda}\right)
\left(\frac{\l}{\mu^2}\right)^\frac14 (|I| \lambda)^\frac14 (|I|\mu^2)^{-\frac14} 
\]
which is favorable if $|I| \geq \mu^{-2}$.

b) $S_\lambda  A \nabla S_\mu u S_\lambda \bar v$ with $\mu \ll
\lambda$. Then using \eqref{lmla} yields a constant 
\[
\mu^\frac32  \lambda^{-2} \ln \lambda \leq \lambda^{-\frac12} \ln \lambda
\]
and a power of $|I|$ can be easily gained as in Remark~\ref{tdelta}.

c) $S_\nu  A \nabla S_\nu u S_\lambda \bar v$ with $\nu \gg
\lambda$. Then we can use \eqref{lmla} but only on $\nu^{-1}$ time
intervals. We obtain a constant
\[
\lambda^\frac12 \nu^{-1} \ln \nu 
\]
and again a power of $I$ is gained as in Remark~\ref{tdelta}.

For the second term in $f_\lambda$ by duality we
need to prove that
\[
\left| \int_I \int_{\R^3}  A^2  u   S_\lambda \bar v dx
dt \right| \lessa   (\lambda |I|)^\delta  \|v\|_{V^2_{A,\sqrt\l} L^2}
\sup_{\lambda}  \| S_\lambda u\|_{U^2_{A,\sqrt\l}(I;L^2)}
\]
We consider the corresponding dyadic pieces
\[
\begin{split}
 \left| \int_I \int_{\R^3}  S_{\lambda_1}A S_{\lambda_2}  A S_{\lambda_3} u   S_\lambda \bar v dx
dt \right| \!\!\!\!& \\ 
 \leq (\lambda  |I|)^{\frac{5}{12}}      \|  &S_{\lambda_1}A\|_{L^6 L^3}
\|  S_{\lambda_2}A\|_{L^6 L^3}         
\| S_{\lambda_3} u\|_{L^\infty L^2} 
\| S_\lambda  v\|_{L^4 L^3} 
\\ \lessa   (\lambda |I|)^{\frac{5}{12}}  &  \lambda_1^{-\frac23}
\lambda_2^{-\frac23}  \|  A\|_{U^2_WH^1}^2
\| S_{\lambda_3} u\|_{L^\infty L^2} 
\|v\|_{V^2_{A,\sqrt\l} L^2} 
\end{split}
\]
where for $v$ we have used the short time Strichartz estimates.
The summation with respect to $\lambda_1$ and $\lambda_2$ is trivial.
So is the summation with respect to $\lambda_3$ since the integral is
zero unless either $\lambda_3=\lambda$ or $\lambda_3 \leq \max\{\lambda_1,\lambda_2\}$.

Finally for the commutator term in $f_\lambda$ we have the bound
\[
\begin{split}
\| [S_\lambda , A_{< \sqrt \lambda} ] \nabla \tilde S_\l u\|_{L^1 L^2}
&\ \lesssim
\l^{-1} (\l |I|) \| [S_\lambda , A_{< \sqrt \lambda} ] \nabla \tilde S_\l u\|_{L^\infty L^2}
\\
&\ \lesssim  \l^{-1} (\l |I|) \| \nabla A_{< \sqrt \lambda} \|_{L^\infty} \|\tilde
S_\lambda u\|_{L^\infty L^2} 
\\
&\ \lesssim  \l^{-\frac14} (\l |I|) \|A\|_{L^\infty H^1} \|\tilde
S_\lambda u\|_{L^\infty L^2}
\end{split}
\]

which again suffices by duality.

b) By the same argument as in part (a) we obtain
\[
\| A_{<\sqrt \lambda} \nabla u - A_{<\nu} \nabla u
\|_{DU^2_{A,\sqrt\lambda}(I;L^2)} \lessa (\lambda |I|)^\delta \| u\|_{U^2 _{A,\sqrt\lambda}(I; L^2)}
\]
which shows that  the two flows are close. Then
\eqref{slec} follows due to Lemma~\ref{flowdiff}.
\end{proof}

\section{ Generalized wave packet decompositions
and long range trilinear estimates}
\label{largescale}

Denote by $C_1(\mu,\lambda,|I|)$ 
the best constant in the estimate
\[
\left| \int_{I} \int_{\R^3} B_\mu     u_\l   \bar{v}_\l dx
dt \right| \leq  C_1(\mu,\lambda,|I|) \| B_\mu \|_{U^2_W L^2}  \|u_\l
\|_{U^2_{A,\l,\l} L^2} \| v_\l \|_{U^2_{A,\l,\l} L^2}
\]
with $B_\mu$, $u_\l$ and $v_\l$ localized at frequencies $\mu$, $\l$,
respectively $\l$. Similarly, let $C_2(\mu,\lambda,|I|)$ be the best
constant in the estimate
\[
\left| \int_{I} \int_{\R^3}  B_\l   u_\mu   \bar{v}_\l dx
dt \right| \leq  C_2(\mu,\lambda,|I|) \| B_\l\|_{U^2_W L^2}  \|u_\mu
\|_{U^2_{A,\mu,\mu} L^2} \| w_\l\|_{U^2_{A,\l,\l} L^2}
\]
with $B_\l$, $u_\mu$ and $v_\l$ localized at frequencies $\l$, $\mu$,
respectively $\l$. 

As a consequence of Proposition~\ref{u2tri} we have 
\[
C_1(\mu,\lambda,\lambda^{-1}) \lessa   \lambda^{-\frac12}
 \min{\{\mu \lambda^{-\frac12},1\}}
\]
A trivial summation shows that for larger time intervals we have
\begin{equation}
C_1(\mu,\lambda,|I|) \lessa   \lambda^{-\frac12}
\min\{\mu \lambda^{-\frac12},1\}
(1+\lambda |I|)
\label{tri0}\end{equation}
We seek to iteratively improve this  to 
\begin{equation} \label{righttri}
C_1(\mu,\lambda,|I|) \lessa    \lambda^{-\frac12}
\min{\{\mu \lambda^{-\frac12},1\}}(1+\lambda |I|)^\frac12
\end{equation}
for intervals $I$ almost up to length $1$. To achieve this we
iteratively produce an increasing sequence of times $T$ for which
\eqref{righttri} holds for $|I| \leq T$. At the same time we seek to
improve $C_2(\mu,\lambda,|I|)$ in a similar manner, as well as extend
the time for which the local energy and local Strichartz estimates
hold.  More precisely, we consider a set of properties as follows:

(i) Paradifferential approximation of the flow:
\begin{equation}
\|S_\lambda  u\|_{U^2_{A,\l,\l}(I; L^2)} \lessa \|u\|_{U^2_{A}(I, L^2)}
\label{rt0}\end{equation}

(ii) Trilinear bounds:
\begin{equation}
C_1(\mu,\lambda,|I|) \lessa  \lambda^{-\frac12}
\min{\{\mu \lambda^{-\frac12},1\}}(1+\lambda |I|)^\frac12
\label{rt1}\end{equation}
\begin{equation}
C_2(\mu,\lambda,|I|) \lessa \  \mu^\frac12 
\lambda^{-1} (1+\lambda |I|)^\frac12
\label{rt2}\end{equation}

(iii) Local energy and local Strichartz estimates for each cube $Q$
of size $1$ and each function $u_\l$ localized at frequency $\lambda$:
\begin{equation}
\| u_\l\|_{L^2 (I;L^2(Q))} \lessa \lambda^{-\frac12}
 \|u_\l\|_{U^2_{A,\l,\l} L^2}, 
\label{rt3}\end{equation}
\begin{equation}
\| u_\l\|_{L^2 (I;L^6(Q))} \lessa \|u_\l\|_{U^2_{A,\l,\l} L^2}.
\label{rt4}\end{equation}

So far, by Propositions~\ref{pst},\ref{paal},\ref{u2tri} we know that 
the above estimates ~\eqref{rt0}-\eqref{rt4} hold if $|I| \leq
\lambda^{-1}$. On the other hand, in order to prove Theorem~\ref{ta}
we need to know that  ~\eqref{rt0}-\eqref{rt4} hold if $|I| \leq
\lambda^{-\epsilon}$ for $\e$ arbitrarily small (see also Remark~\ref{tdelta}).
 This is accomplished in the next result.

\begin{p1} 
  Let $0 < \alpha \leq 1$.  Assume that the estimates
  \eqref{rt0}-\eqref{rt4} hold for $|I| < \lambda^{-\alpha}$.  Then
  \eqref{rt0}-\eqref{rt4} hold for $|I| < \lambda^{-\beta}$ for each
  $\beta > \frac34 \alpha$.
\end{p1}

\begin{proof}
We first improve the time range of the paradifferential calculus:
\begin{l1}
a)  For each frequency $\lambda$ we have
\[
\| S_\lambda u\|_{U^2_{A,\l,\l}(I, L^2)} \lesssim_A \| u\|_{U^2_A(I,
  L^2)}, \qquad |I| \leq  \lambda^{ -\frac{\alpha}2}  (\log \lambda)^{-\frac32}
\]

b) For each frequency $ \lambda^{1-\frac{\alpha}2}  \log \l < \nu  \ll  \lambda$ we have
\[
\| u\|_{U^2_{A,\nu,\lambda}(I, L^2)} \approx_A  \| u\|_{U^2_{A,\l,\l}(I, L^2)},
\qquad |I| \leq T(\lambda,\nu)=  \nu \lambda^{-1 -\frac{\alpha}2}  (\log \lambda)^{-1}
\]

\end{l1}
\begin{proof}
a) We observe that \eqref{rt1}, \eqref{rt2}, \eqref{rt3}  for $|I| =
\lambda^{-\alpha}$ trivially lead to bounds for longer time,
\begin{equation}
C_1(\mu,\lambda,|I|) \lessa   \lambda^{-\frac12}
\min{\{\mu \lambda^{-\frac12},1\}}
(1+\lambda |I|)^\frac12 (1+ \lambda^\alpha |I|)^\frac12
\label{rt1a}\end{equation}
\begin{equation}
C_2(\mu,\lambda,|I|) \lessa  \mu^\frac12 \lambda^{-1} 
(1+\lambda |I|)^\frac12 (1+ \lambda^\alpha |I|)^\frac12
\label{rt2a}\end{equation}
\begin{equation}
\| u_\l\|_{L^2 (I,L^6(Q))} \lessa
 (1+ \lambda^\alpha |I|)^\frac12 \|u_\l\|_{U^2_{A,\l,\l} L^2}
\label{rt4a}\end{equation}
Furthermore, due to \eqref{rt0}, we obtain the same  constant 
in the trilinear estimates when we replace $u_\l$ by $S_\l u$
and $\| u_\l\|_{U^2_{A,\l,\l} L^2}$ by $\|u\|_{U^2_A L^2}$. Also 
by the argument in Lemma~\ref{principle}, we can also replace one of
the $U^2$ norms with a $V^2$ norm at the expense of an additional
$\ln \lambda$ loss.

The rest of the proof is similar to the proof of Proposition~\ref{plla}. 
For $u$ solving \eqref{schra} we write
\[
\left(i \partial_ t - \Delta + i A_{\ll \lambda} \nabla \tilde
  S_\lambda + i \tilde S_\lambda A_{\ll \lambda} \nabla \right)
S_\lambda u = f_\lambda
\]
where 
\[
f_\lambda = S_\lambda ( 2i A_{\gtrsim \lambda} \nabla u + A^2 u) + i [
S_\lambda , A_{\ll \lambda} ] \tilde
  S_\lambda \nabla u 
\]
Hence it suffices to prove that
\[
\| f_\lambda\|_{DU^2_{A,\l,\l}(I;L^2)} \lessa \| u\|_{ U^2_{A}(I,L^2)}
\]
We use duality and consider each term in $f_\lambda$. For the first
one we need to show that
\[
\left| \int_{I} \int_{\R^3}  A_{\gtrsim \lambda}  \nabla  u  S_\lambda \bar{v} dx
dt \right| \lessa \|A\|_{U^2_W H^1}  \| u\|_{ U^2_{A}(I,L^2)} \|
v\|_{V^2_{A,\l,\l} L^2}
\]
After a Littlewood-Paley decomposition of the first two factors we
need  to consider the following three cases for the integrand:

(i) $S_\lambda A \nabla S_\lambda u   S_\lambda \bar{v}$. Then we use 
\eqref{rt1a} to obtain a constant
\[
\log \lambda \ \lambda^\frac12 (\lambda |I|)^\frac12 (\lambda^\alpha
|I|)^\frac12  \lambda^{-1} = \log{\lambda} \ \lambda^{\frac{\alpha}2} |I| 
\]
which is satisfactory given the range for $I$.

(ii) $S_\lambda A \nabla S_\mu u   S_\lambda \bar{v}$, $\mu \ll
\lambda$. Then we use 
\eqref{rt2a} to obtain a constant
\[
  \mu^\frac32 \lambda^{-1} (\lambda |I|)^\frac12 (\lambda^\alpha
|I|)^\frac12  \lambda^{-1} \leq \lambda^{\frac{\alpha}2} |I| 
\]
which is much better than we need.

(iii) $S_\nu A \nabla S_\nu u S_\lambda \bar{v}$, $\lambda \ll \nu$.
Then we use \eqref{rt2a} to obtain a constant
\[
  \lambda^\frac32 \nu^{-1} (\nu |I|)^\frac12 (\nu^\alpha
|I|)^\frac12  \nu^{-1} \leq \lambda^{\frac{\alpha}2} |I| 
\]

For the second term in $f_\lambda$ by duality we need to prove that
\[
\left| \int_{I} \int_{\R^3}  A B   u  S_\lambda \bar{v} dx
dt \right| \lessa \|A\|_{U^2_W H^1} \|B\|_{U^2_W H^1} \| u\|_{ U^2_{A}(I,L^2)} \|
v\|_{V^2_{A,\l,\l} L^2}
\]
Due to the finite speed of propagation for the wave equation, see
\eqref{wcloc}, it suffices to consider the case when $A$ and $B$ are
supported in a unit cube. Then we bound $A$ and $B$ in $L^\infty L^6$,
$u$ in $L^2 L^6$ as in \eqref{rt4a}, and $v$ in $L^\infty L^2$ and use
Holder's inequality with respect to time.  This yields the same constant
$\lambda^{\frac{\alpha}2} |I|$ as above.

Finally we consider the commutator term in  $f_\lambda$. This can be
represented in the form of a rapidly convergent series of the form
\[
[S_\lambda,A_{\ll \lambda}] \nabla \tilde S_\lambda u =  \sum_j
S_\lambda^{1j} (\nabla A_{\ll \lambda} S_\lambda^{2j} u)
\]
where $S^{1j}_\l, S^{2j}_\l$ are operators similar to $S_\l$.

Then by duality it suffices to prove that
\[
\left| \int_I\! \int_{\R^2}\!\! \nabla A_{\ll \lambda} \nabla S_\l u
  S_\lambda \bar{v} dx dt \right| \! \lessa \!  (\log{\lambda})^\frac32
\l^{\frac{\alpha}{2}} |I|  
\| \nabla A  \|_{U^2_W L^2} \| u\|_{ U^2_{A} L^2} \| v\|_{V^2_{A,\l,\l} L^2}
\]
For this it suffices to consider a Littlewood-Paley decomposition of
$A_{\ll \lambda}$ and to apply \eqref{rt1a} for each dyadic piece.

% \[
% \l^{-1} \sum_{\mu \ll \l} \left| \int_I \int_{\R^2} \nabla A_{\mu} \nabla \tilde S_\l u  \bar{v} dx dt \right| 
% \]

% \[
% \les  \l^{-1} \log{\lambda} (\l |I|)^\q (\l^\a |I| )^\q \sum_{\mu \ll \l} \min{(\mu,\l^\q)} \| \nabla A_\mu  \|_{U^2_W L^2} \| u\|_{ U^2_{A}(I,L^2)} \| v\|_{V^2_{A,\l,\l} L^2}
% \]

% \[
% \les  (\log{\lambda})^\frac32  \l^{\frac{\a}{2}} |I|  \| \nabla A  \|_{U^2_W L^2} \| u\|_{ U^2_{A}(I,L^2)} \| v\|_{V^2_{A,\l,\l} L^2}
% \]

b) By virtue of Lemma~\ref{flowdiff} it suffices to show that
\[
\| A_{\ll \lambda}  \nabla \tilde S_\lambda u - A_{<\nu} \nabla \tilde
S_\lambda u\|_{DU^2_{A,\l,\l}} \lessa \|u\|_{U^2_{A,\l,\l}}
\]
and the similar bound for $ \tilde S_\lambda A_{\ll \lambda} \nabla
-\tilde S_\lambda A_{<\nu} \nabla$.  By duality this becomes
\[
\left| \int_{I} \int_{\R^3} \sum_{\mu = \nu}^\lambda S_\mu A  \nabla \tilde S_\lambda u
\bar{v} dx dt \right| \lessa \| A\|_{U^2_A H^1}
\|u\|_{U^2_{A,\l,\l}} \|v\|_{V^2_{A,\l,\l}}
\]
Indeed, for $|I| > \lambda^\alpha$ the estimate \eqref{rt1a} yields a constant
\[
\log \lambda \ \lambda^\frac12   (\lambda |I|)^\frac12 (\lambda^\alpha
|I|)^\frac12 \nu^{-1} = \log \lambda \ \nu^{-1} \lambda^{\frac{\alpha}2+1} |I|
\]
which leads to the restriction 
\[
|I| \leq T(\lambda,\nu)= (\log \lambda)^{-1} \nu \lambda^{-\frac{\alpha}2-1}
\]
We observe that this is useful only if it provides information on 
time intervals with $|I| > \lambda^{-\alpha}$. This leads to the
condition
\[
\nu >  \l^{1-\frac{\alpha}2} \log \l.
\]

\end{proof}
 
Next we consider the \eqref{schranu} evolution and we construct a
generalized wave packet structure for the flow. The frequency scale is
$\delta \xi = \nu$ and the time scale is $T(\lambda,\nu)$ therefore it
is natural to define the spatial scale by $\delta x = \nu
T(\lambda,\nu)$, as in the case of the flat flow.

We first partition the initial data. Let $\phi$ be a smooth unit bump
function in $\R^3 \times \R^3$ so that 
\[
\sum_{k,j \in \Z^n}  \phi(x-k,\xi-j) = 1
\]
 Denote
\[
\phi_{x_0,\xi_0}^\nu (x,\xi) =  \phi\left(\frac{x-x_0}{\nu T(\lambda,\nu)},\frac{\xi-\xi_0} \nu\right)
\]
where 
\[
(x_0,\xi_0) \in \Z_\nu^{2 \times 3} = (\nu T(\lambda,\nu) \Z)^3 \times (\nu
\Z)^3
\]
Then consider an almost orthogonal decomposition of the initial 
data
\[
u_0 = \sum_{(x_0,\xi_0) \in \Z_\nu^{2 \times 3}} \phi_{x_0,\xi_0}^\nu (x,D) u_0
\]
We denote the corresponding solutions to \eqref{schranu} by
$u_{x_0,\xi_0}$, and we call them generalized wave packets.
 To measure their evolution we consider the family of
operators
\[
L_{x_0,\xi_0} = \left\{ \nu^{-1} (\xi-\xi_0),\ (\nu
  T(\lambda,\nu))^{-1} (x - x_0 - 2 t\xi)  \right\}
\]
which commute with $i \partial_t -\Delta$. The following lemma
shows that $u_{x_0,\xi_0}$ is concentrated in a tube
\[
T^\nu_{x_0,\xi_0} =  \{ (x,\xi): |x-x_0-2t\xi_0| \leq \nu T(\l, \nu),
|\xi-\xi_0| \leq \nu \}
\]
and decays rapidly away from it.

\begin{l1} \label{wpnu}
The solutions $u_{x_0,\xi_0}$ for the \eqref{schranu} flow satisfy
\begin{equation}
  \sum_{(x_0,\xi_0) \in \Z_\nu^{2 \times 3}} \sum_{|\alpha| \leq N}
  \|L_{x_0,\xi_0}^\alpha u_{x_0,\xi_0}(t)\|_{U^2_{A,\nu,\lambda}
    (I,L^2)}^2 
\lessa  \|u_0\|_{L^2}^2, \qquad |I| \leq T(\lambda,\nu)
\label{genwp}\end{equation}
\end{l1}
\begin{proof}
At time $0$ we clearly have
\[
 \sum_{|\alpha| \leq N}
\|L_{x_0,\xi_0}^\alpha u_{x_0,\xi_0}(0)\|_{L^2}^2 
\lesssim \|u_0\|_{L^2}^2
\]
therefore it suffices to prove that a single generalized wave packet
satisfies
\begin{equation}
 \sum_{|\alpha| \leq N}
\|L_{x_0,\xi_0}^\alpha u_{x_0,\xi_0}\|_{U^2_{A,\nu} (I,L^2)}^2 
\lessa \sum_{|\alpha| \leq N}
\|L_{x_0,\xi_0}^\alpha u_{x_0,\xi_0}(0)\|_{L^2}^2 
\label{genwpa}\end{equation}
This  follows iteratively from 
\begin{equation}
\|L_{x_0,\xi_0} v\|_{U^2_{A,\nu,\lambda} (I,L^2)}^2 
\lessa 
\|L_{x_0,\xi_0} v (0)\|_{L^2}^2 + \| v\|_{U^2_{A,\nu.\lambda} (I,L^2)}^2
\label{genwpb}\end{equation}
for which, in turn, we need the commutator bound
\begin{equation}
\|[ L_{x_0,\xi_0}, A_{<\nu} \tilde S_{\lambda} \nabla] v 
\|_{DU^2_{A,\nu,\lambda} (I,L^2)}
\lesssim_A \| v\|_{U^2_{A,\nu,\lambda} (I,L^2)}^2
\label{genwpc}\end{equation}
as well as the similar one for the operator $\tilde
S_{\lambda}A_{<\nu} \nabla$.

If $L_{x_0,\xi_0} = \nu^{-1} (\xi-\xi_0)$ then 
\[
[ L_{x_0,\xi_0}, A_{<\nu} \tilde S_{\lambda} \nabla] =  \nu^{-1} (\nabla A_{<\nu}) \tilde S_{\lambda} \nabla
\]
therefore by duality we need to show that
\[
\left| \int_{I} \int_{\R^3}   \nabla A_{<\nu}  \nabla \tilde S_\lambda u
\bar{v} dx dt \right| \lessa \nu \| A\|_{U^2_A H^1}
\|u\|_{U^2_{A,\nu,\lambda}} \|v\|_{V^2_{A,\nu,\lambda}}
\]
which follows from \eqref{rt1a}.

If $L_{x_0,\xi_0} =  (\nu
  T(\lambda,\nu))^{-1} (x - x_0 - 2 t\xi)$ then 
\[
[ L_{x_0,\xi_0}, A_{<\nu} \tilde S_{\lambda} \nabla] =  (\nu
  T(\lambda,\nu))^{-1} (A_{<\nu} + t (\nabla A_{<\nu})\nabla) \tilde S_{\lambda} 
\]
The second term is as above. For the first by duality we need to show that
\[
\left| \int_{I} \int_{\R^3}    A_{<\nu}   \tilde S_\lambda u
\bar{v} dx dt \right| \lessa \nu T(\lambda,\nu)   \| A\|_{U^2_A H^1}
\|u\|_{U^2_{A,\nu,\lambda}} \|v\|_{V^2_{A,\nu,\lambda}}
\]
which is much weaker and follows again from \eqref{rt1a}.

\end{proof}

The parameter $\nu$ is chosen so that the packets move away from their
initial support by the time $\l^{-\alpha}$. For this we impose the condition
\[
\nu T(\lambda,\nu) < \lambda^{1-\alpha -\epsilon}
\]
for some $\epsilon > 0$. This is satisfied if we choose $\nu$ of the
form 
\[
\nu = \lambda^{1-\frac{\alpha}4 -\epsilon}
\]
in which case we have 
\[
\nu T(\lambda,\nu) = \lambda^{-\frac{3\alpha }{4} -2\e} (\log \l)^{-1}
\]
This leads to the choice of $\beta$ in the proposition. Indeed, the proof 
of the proposition is concluded due to the following lemma: 

\begin{l1} 
  Let $\e > 0$.  Choose $\nu$ so that $\nu T(\lambda,\nu) <
  \lambda^{1-\alpha-\e}$ and $\nu < \lambda^{1-\e}$.  Then
  \eqref{rt1}-\eqref{rt4} hold for $|I| < T(\lambda,\nu)$.
\end{l1}
\begin{proof}
  By the previous lemma we can replace the $(A_{\ll \lambda},\lambda)$
  flow by the $(A_{<\nu},\lambda)$ flow in \eqref{rt1}-\eqref{rt4}.
  We begin with \eqref{rt1}. Without any restriction
  in generality we can assume that $u$ and $v$ are
  $U^2_{A,\nu,\lambda} L^2$ atoms. The generalized wave packet
  decomposition in Lemma~\ref{wpnu} easily extends to
  $U^2_{A,\nu,\lambda} L^2$ atoms. Indeed, if we denote 
by $S_\nu(t,s)$ the evolution generated by \eqref{schranu} then an atom $u_\l$
of the form
\[
u_\l(t) = \sum_k 1_{[t_k,t_{k+1})} S_\nu(t,t_k) u^k_\l 
\]
can be partitioned as
\begin{equation}
u = \sum_{x_0} \sum_{|\xi_0| \approx \l}  u_{x_0,\xi_0}
\label{u=u}\end{equation}
where 
\begin{equation}
 u_{x_0,\xi_0}(t) = \sum_k 1_{[t_k,t_{k+1})} S_\nu(t,t_k)
 u^k_{x_0,\xi_0}
\label{u=uk}\end{equation}
with
\[
 u^k_{x_0,\xi_0} =  \phi^\nu_{x_0+2t_k\xi_0,\xi_0} (x,D) u^k_\l
\]
Then, denoting
\[
\tri u_{x_0,\xi_0} \tri^2 = \sum_{k} \sum_{|\alpha|
  \leq N}   \|L_{x_0,\xi_0}^\alpha (t_k) u^k_{x_0,\xi_0}\|_{L^2}^2 
\]
 we have the orthogonality relation
\begin{equation}
\sum_{x_0} \sum_{|\xi_0| \approx \l}  \tri u_{x_0,\xi_0} \tri^2
\lesssim \sum_k \|u_k\|_{L^2}^2
\label{u=uka}\end{equation}
We argue in a similar manner for $v$. 
Hence it suffices to take $u_\l$  as in \eqref{u=u},\eqref{u=uk}, and 
similarly for $v_\l$, and prove that for $|I| \leq T(\lambda,\nu)$ we have
\begin{equation}
\begin{split} \left| \int_{I} \int_{\R^3} B_\mu     u_\l   \bar{v}_\l dx
dt \right| &\  \lessa RHS\eqref{rt1} \cdot 
\| B_\mu \|_{U^2_W L^2} 
\\ &\ \!\!\!\!\! \left(\sum_{x_0}  \sum_{|\xi_0| \approx \l} \tri
u_{x_0,\xi_0}\tri^2\right)^\frac12
\left(\sum_{x_0}  \sum_{|\xi_0| \approx \l} \tri
v_{x_0,\xi_0}\tri^2\right)^\frac12
\end{split}
\label{triest}\end{equation}
We proceed with several reductions, which will eventually lead to
shorter time intervals.

{\bf 1. Reduction to spatial scale $\l T(\lambda,\nu)$.} Heuristically,
in time $T(\lambda,\nu)$ the frequency $\lambda$ waves for
\eqref{schranu} travel by $\lambda T(\lambda,\nu)$.  Hence we
partition the space into cubes $\{Q_j\}_{j \in \Z^3}$ of size $\lambda
T(\lambda,\nu)$. Correspondingly, we decompose $u_\l$ into
\[
u_\l = \sum_{j \in \Z^3} u_j, \qquad u_j =  \sum_{x_0 \in Q_j}  u_{x_0,\xi_0}
\]
and similarly for $v$.  
By \eqref{genwp} the functions $u_j$ decay rapidly away
from an enlargement of $Q_j$. Precisely, if $x_0 \in Q_j$ and 
$|j-k| \geq 10$ then   the separation between the tube
$T_{x_0,\xi_0}^\nu$ and $Q_k$ is $O(|j-k|\lambda T(\lambda,\nu))$. 
Comparing this with the tube thickness $\nu  T(\lambda,\nu))$
 we have
\[
|u_j(t,x)| \lesssim \lambda^{-N} |j-k|^{-N} \sum_{x_0 \in Q_j}
\sum_{|\xi_0| \approx \l} \tri u_{x_0,\xi_0} \tri^2, \quad x \in Q_k,
\ |j-k| \geq 10
\]
Thus in \eqref{triest} it suffices to consider the output of $u_j$ and
$v_{j_1}$ for $|j-j_1| < 20$, and only within an enlargement $C Q_j$
of $Q_j$; the rest is trivially estimated using the above bound.
Furthermore, by Cauchy-Schwartz it suffices to consider a fixed $j$
and $k$.  By a slight abuse of notation we set $k=j$ in the
sequel. Then  \eqref{triest} is reduced to 
\begin{equation}
\begin{split} \left| \int_{I} \int_{\R^3}   \chi_{C Q_j} B_\mu     u_j   \bar{v}_{j} dx
dt \right| &\  \lessa RHS\eqref{rt1} \cdot 
\| B_\mu \|_{U^2_W L^2} 
\\ &\hspace{-1in} \left(\sum_{x_0 \in Q_j}  \sum_{|\xi_0| \approx \l} \tri
u_{x_0,\xi_0}\tri^2\right)^\frac12
\left(\sum_{x_0\in Q_j}  \sum_{|\xi_0| \approx \l} \tri
v_{x_0,\xi_0}\tri^2\right)^\frac12
\end{split}
\label{triesta}\end{equation}

{\bf 2. Reduction to small angles.} Here we partition the $\lambda$
annulus $A_\lambda$ in frequency into small angles of size
$\frac1{10}$ with centers in $\Theta \subset \S^2$,
\[
A_\lambda = \bigcup_{\theta \in \Theta} A_{\lambda,\theta}
\]
Then we divide
\[
u_j = \sum_{\theta \in \Theta} u_{j,\theta}, \qquad u_{j,\theta} = 
\sum_{x_0 \in Q_j}  \sum_{\xi \in A_\theta} u_{x_0,\xi_0}
\]
and similarly for $v_j$. It remains to prove that
\begin{equation}
\begin{split} \left| \int_{I} \int_{\R^3}   \chi_{C Q_j} B_\mu    
 u_{j,\theta}   \bar{v}_{j,\omega} dx
dt \right| &\  \lessa RHS\eqref{rt1} \cdot 
\| B_\mu \|_{U^2_W L^2} 
\\ &\hspace{-1in} \left(\sum_{x_0 \in Q_j}  \sum_{\xi_0 \in A_{\l,\theta}} \tri
u_{x_0,\xi_0}\tri^2\right)^\frac12
\left(\sum_{x_0\in Q_j}  \sum_{\xi_0 \in A_{\l,\omega}}  \tri
v_{x_0,\xi_0}\tri^2\right)^\frac12
\end{split}
\label{triestb}\end{equation}

{\bf 3. Reduction to a spatial strip of size $\lambda^\e \nu
  T(\lambda,\nu)$.} Given directions $\theta$ and $\omega$ as above we
choose a coordinate, say $\xi_1$, so that both the $\theta$ and the
$\omega$ sectors $A_\theta$, $A_\omega$ are away from $\xi_1=0$.
Dividing the spatial coordinates $x = (x_1,x')$ we partition the space
into strips $S_k$ of thickness $\lambda^\e \nu T(\lambda,\nu)$ in the $x_1$
direction.

Arguing as in \eqref{wcloc}, the $A$ factor is square summable with
respect to the strips $S_k$. There are about $\lambda^{1-\epsilon}
\nu^{-1}$ such strips which intersect $30 Q_j$.  Hence by losing a
$\lambda^{\frac{1-\e}2} \mu^{-\frac12}$ factor we can use Holder's
inequality to reduce the problem to the case when $A$ is supported in
a single spatial strip $S_k$. It remains to prove that
\begin{equation}
\begin{split} \left| \int_{I} \int_{\R^3}   \chi_{S_k \cap C Q_j} B_\mu    
 u_{j,\theta}   \bar{v}_{j,\omega} dx
dt \right| &\  \lessa RHS\eqref{rt1} \cdot \lambda^{-\frac{1-\e}2} 
\mu^{\frac12}
\| B_\mu \|_{U^2_W L^2} 
\\ &\hspace{-1.2in} \left(\sum_{x_0 \in Q_j}  \sum_{\xi_0 \in A_{\l,\theta}} \tri
u_{x_0,\xi_0}\tri^2\right)^\frac12
\left(\sum_{x_0\in Q_j}  \sum_{\xi_0 \in A_{\l,\omega}}  \tri
v_{x_0,\xi_0}\tri^2\right)^\frac12
\end{split}
\label{triestc}\end{equation}

{\bf 4. Reduction to spatial scale $\lambda^\e \nu
  T(\lambda,\nu)$.} Due to the choice of coordinates above,
the packets in $u_{j,\alpha}$ and $v_{j,\omega}$ travel 
in directions which are transversal to $S_k$. Hence if we partition
$S_k$ into cubes $\tilde Q_l$ of size $\lambda^\e \nu
  T(\lambda,\nu)$, each packet will intersect only finitely many 
cubes. Then we partition further
\[
u_{j,\theta} = \sum_l u_{j,\theta,l}, \qquad u_{j,\theta,l} =
\sum_{(x_0,\xi_0) \in A_{j,\theta,l}} u_{x_0,\xi_0}
\]
where
\[
A_{j,\theta,l} = \{ (x_0,\xi_0); x_0 \in Q_j,\ \xi_0 \in
A_{\l,\theta}, \ \{x_0 + \R\xi_0\} \cap
  \tilde Q_l \not = \emptyset \}
\]
and packets $u_{x_0,\xi_0}$ intersecting more than one cube
are arbitrarily placed in one of the terms.

Arguing as in the first reduction, for $|l -l_1|\gg 1$ the size of
$u_{j,\theta,l}$ in $\tilde Q_{l_1}$ is rapidly decreasing,
\[
|u_{j,\theta,l}(t,x)|^2 \lesssim \lambda^{-N} |l-l_1|^{-N} 
\sum_{(x_0,\xi_0) \in A_{j,\theta,l}} \tri
u_{x_0,\xi_0}\tri^2
 , \quad x \in \tilde Q_{l_1}, \ |l-l_1| \gg 1
\]
Hence the joint contribution of $u_{j,\theta,l}$ and
$v_{j,\omega,l_1}$ in \eqref{triestc} is nontrivial only on the diagonal
$|l-l_1| \lesssim 1$. By Cauchy-Schwartz with respect to  $l$ it
suffices to estimate this contribution for fixed $l,l_1$, in an
enlarged cube $C \tilde Q_l$.  By a slight abuse of notation we set $l
= l_1$. Then we need to show that
\begin{equation}
\begin{split} \left| \int_{I} \int_{\R^3}   \chi_{C Q_l} B_\mu    
 u_{j,\theta,l}   \bar{v}_{j,\omega,l} dx
dt \right| &\  \lessa RHS\eqref{rt1} \cdot \lambda^{-\frac{1-\e}2} 
\mu^{\frac12}
\| B_\mu \|_{U^2_W L^2} 
\\ &\hspace{-1.2in} \left(  \sum_{(x_0,\xi_0) \in A_{j,\theta,l}} \tri
u_{x_0,\xi_0}\tri^2\right)^\frac12
\left(  \sum_{(x_0,\xi_0) \in A_{j,\omega,l}}  \tri
v_{x_0,\xi_0}\tri^2\right)^\frac12
\end{split}
\label{triestd}\end{equation}

{\bf 5. Reduction to a time interval of size $\lambda^{\e-1} \nu
  T(\lambda,\nu)$.}  Each tube $T_{x_0,\xi_0}$ intersects the cube
$Q_l$ in a time interval of size $\lambda^{\e-1} \nu T(\lambda,\nu)$.
Hence it is natural to partition the time interval $I$ into
subintervals $I_m$ of length $\lambda^{\e-1} \nu
  T(\lambda,\nu)$. Correspondingly, we split $u_{j,\theta,l}$ into
\[
u_{j,\theta,l} =\sum_m u_{j,\theta,l,m}, \qquad u_{j,\theta,l,m} = 
\sum_{(x_0,\xi_0) \in A_{j,\theta,l,m}} u_{x_0,\xi_0}
\]
where
\[
A_{j,\theta,l,m} = \{ (x_0,\xi_0); x_0 \in Q_j,\ \xi_0 \in
A_{\l,\theta}, \ T^\nu_{x_0,\xi_0} \cap
  I_m \times  \tilde Q_l \not = \emptyset \}
\]
and similarly for $v_{j,\omega,l}$. Again, the size of
$u_{j,\theta,l,m}$ in $I_{m_1} \times Q_l$ is negligible if $|m-m_1|
\gg 1$. Thus by Cauchy-Schwartz with respect to $m$ the estimate 
 \eqref{triestd} reduces to the case of a single interval $I_m$,
\begin{equation}
\begin{split} \!\!\!\left| \int_{I_m} \! \int_{\R^3}   \chi_{C Q_l} B_\mu    
 u_{j,\theta,l,m}   \bar{v}_{j,\omega,l,m} dx
dt \right| &  \lessa\! RHS\eqref{rt1}  \lambda^{-\frac{1-\e}2} 
\mu^{\frac12}
\| B_\mu \|_{U^2_W L^2} 
\\ &\hspace{-1.4in} \left(  \sum_{(x_0,\xi_0) \in A_{j,\theta,l,m}} \tri
u_{x_0,\xi_0}\tri^2\right)^\frac12
\left(  \sum_{(x_0,\xi_0) \in A_{j,\omega,l,m}}  \tri
v_{x_0,\xi_0}\tri^2\right)^\frac12
\end{split}
\label{trieste}\end{equation}
But this follows from the hypothesis since
\[
|I_m| = \lambda^{\e-1} \nu  T(\lambda,\nu) \leq \lambda^{-\alpha}
\]
and
\[
|I|^\frac12 \lambda^{-\frac{1-\e}2} 
\mu^{\frac12} = |I_m|
\]

In the case of \eqref{rt2} the argument is similar but with several
adjustments which we outline. 

{\bf 1. Reduction to spatial scale $ \max\{\mu,\l^\e \nu\} T(\lambda,\nu)$.}
This smaller initial localization scale is possible since frequency
$\nu$ Schr\"odinger waves travel with speed $\mu$, so within time 
$ T(\lambda,\nu)$ they can spread only as far as $\mu
T(\lambda,\nu)$. Thus for the $u_\mu$ factor we have square 
summability on the $\mu T(\lambda,\nu)$. For the wave factor $B_\l$
by \eqref{wcloc} we have square summability on the same scale,
therefore we are allowed to localize spatially the estimate  on the 
$\mu T(\lambda,\nu)$ scale. If $\mu$ is small we only can take partial
advantage of this due to the wider spread of frequency $\l$ pakets.

{\bf 2. Reduction to small angles.} This is as before.

{\bf 3. Reduction to a spatial strip of size $\lambda^\e \nu
  T(\lambda,\nu)$.} The difference here is that we have only about
$\max\{ 1, \mu \l^{-\e} \nu^{-1}\}$ strips intersecting a 
 $\max\{\mu,\l^\e \nu\} T(\lambda,\nu)$ cube, therefore we only loose
a factor of 
\[
\max\{ 1, \mu^\frac12 \l^{-\frac{\e}2}
\nu^{-\frac12}\}
\]

{\bf 4. Reduction to spatial scale $\lambda^\e \nu
  T(\lambda,\nu)$.} This is as before.

{\bf 5. Reduction to a time interval of size $\lambda^{\e-1} \nu
  T(\lambda,\nu)$.} As before, the frequency $\lambda$ wave packets
spend a time $ \lambda^{\e-1} \nu T(\lambda,\nu)$ inside a $
\lambda^\e \nu T(\lambda,\nu)$ cube $Q$. However, the frequency $\mu$
packets spend a longer time $\mu^{-1}\lambda^{\e} \nu T(\lambda,\nu)$
inside $Q$.  Hence we can carry out first a lossless reduction down to
time scale $\min\{ 1, \mu^{-1}\lambda^{\e} \nu \} T(\lambda,\nu)$. To
further reduce the time scale to $\lambda^{\e-1} \nu T(\lambda,\nu)$
we can only use the square summability for the frequency $\l$ waves,
therefore we apply Cauchy-Schwartz and loose an additional factor
of
\[
\min\{ 1, \mu^{-1}\lambda^{\e} \nu \}^\frac12 (\lambda^{\e-1} \nu)^{-\frac12}
\]

Finally, combining the two losses in Steps 3 and 5 we obtain a total loss of
\[
 (\lambda^{\e-1} \nu)^{-\frac12}
\]
which is identical to the one in Step 3  of the proof of
\eqref{rt1}. We conclude as before.

The argument is considerably  simpler in the case of \eqref{rt3} and
\eqref{rt4}. There each packet intersects $Q \times I$ in a time
interval which is shorter than $ \lambda^{-1} \nu T(\lambda, \nu) <
\lambda^{-\alpha}$. Grouping the wave packets with respect to such
time intervals we obtain the square summability of the outputs
and reduce the problem to the shorter time scale $\lambda^{-\alpha}$.

\end{proof}

\end{proof}

\bibliographystyle{plain}
\bibliography{MS}

\end{document}